\documentclass[12pt]{article}
\usepackage{amsmath,amssymb,enumerate,bm}
\usepackage{graphicx,xcolor}
\usepackage[all]{xy}

\makeatletter
\newcommand{\subscripts}[3]{%
  \@mathmeasure\z@\displaystyle{#2}%
  \global\setbox\@ne\vbox to\ht\z@{}\dp\@ne\dp\z@
  \setbox\tw@\box\@ne
  \@mathmeasure4\displaystyle{\copy\tw@_{#1}}%
  \@mathmeasure6\displaystyle{{#2}_{#3}}%
  \dimen@-\wd6 \advance\dimen@\wd4 \advance\dimen@\wd\z@
  \hbox to\dimen@{}\mathop{\kern-\dimen@\box4\box6}%
}
\makeatother
\newcommand{\qed}{~~\hbox{\rule{4pt}{8pt}}}

\newtheorem{Thm}{Theorem}[section]

\newtheorem{Lem}{Lemma}[section]

\newtheorem{Rem}{Remark}[section]

\newcommand{\vv}{\bm{\mathrm{v}}}
\newcommand{\x}{\bm{\mathrm{x}}}

\newcommand{\h}{\mathrm{H}}

\newcommand{\A}{\mathcal{A}}
\newcommand{\T}{\mathcal{T}}
\newcommand{\nn}{\nonumber}
\newcommand{\Prob}{\bm{\mathrm{P}}}
\newcommand{\Probq}{\bm{\mathrm{Q}}}

\newcommand{\X}{\bm{\mathrm{X}}}
\newcommand{\Y}{\bm{\mathrm{Y}}}
\newcommand{\e}{\mathrm{e}}
\newcommand{\p}{\mathfrak{p}}
\newcommand{\m}{\mathfrak{m}}
\newcommand{\ve}{\varepsilon}

\newcommand{\Ker}{\mathrm{Ker}\,}

\newcommand{\del}{\partial}
\newcommand{\ol}{\overline}
\newcommand{\la}{\langle}
\newcommand{\La}{\langle\!\langle}
\newcommand{\ra}{\rangle}
\newcommand{\Ra}{\rangle\!\rangle}
\newcommand{\LA}{\longrightarrow}

\newcommand{\Hom}{\mathrm{Hom}}

\makeatletter
  
  \@addtoreset{equation}{section}
\makeatother

\begin{document}
\title{A remark on a central limit theorem for non-symmetric random walks
on crystal lattices}
\author{
Ryuya Namba
\thanks{Graduate School of Natural Sciences, Okayama University, 3-1-1, 
Tsushima-Naka, Kita-ku, Okayama 700-8530, Japan 
({\tt rnamba@s.okayama-u.ac.jp})
}
}
\date{\today} 	                                         
\maketitle

\begin{abstract}
Recently, Ishiwata, Kawabi and Kotani \cite{IKK} proved two kinds of 
central limit theorems for non-symmetric random walks on crystal lattices
from the view point of discrete geometric analysis developed by Kotani and Sunada.
In the present paper, we establish yet another kind of 
the central limit theorem for them.
Our argument is based on a  
measure-change technique due to 
Alexopoulos \cite{A}. 

\vspace{1.5mm}
\noindent
{\bf 2010 AMS Classification Numbers:} 60J10, 60F05, 60G50, 60B10.

\vspace{1.5mm}
\end{abstract}


\section{Introduction and results}

Let $X=(V, E)$ be a locally finite, connected and oriented graph. 
Here $V$ is the set of all vertices and $E$ the set of all oriented edges. 
For an edge $e \in E$, we denote by $o(e), t(e)$ and $\ol{e}$ the origin, the terminus
and the inverse edge of $e$, respectively. 
We denote by $E_x$ the collection of all edges whose origin is $x \in V$. 
A path $c$ in $X$ with length $n$ is a sequence $c=(e_1, \dots, e_n)$
of edges $e_i$ with $t(e_i)=o(e_{i+1}) \, (i=1, \dots, n-1)$.
We denote by $\Omega_{x, n}(X) \, (x \in V, \, n \in \mathbb{N} \cup \{\infty\})$
the set of all paths of length $n$ for which origin $o(c)=x$.
We also denote by $o(c), t(c)$ 
the origin and the terminus of the path $c$. 
For simplicity, we write $\Omega_x(X):=\Omega_{x, \infty}(X)$.

Let $p : E \LA (0, 1]$ be a {\it transition probability}, 
that is, a positive function on $E$ satisfying 
\begin{equation}\label{trans-prob}
\sum_{e \in E_x}p(e)=1 \qquad (x \in V).
\end{equation}
This induces the probability measure $\mathbb{P}_x$
on $\Omega_x(X)$ and the random walk associated with the transition probability
$p$ is the time homogeneous Markov chain 
$(\Omega_x(X), \mathbb{P}_x, \{w_n\}_{n=0}^\infty)$
with values in $X$ defined by 
$w_n(c):=o(e_{n+1})$ $\big(c=(e_1, e_2, \dots) \in \Omega_x(X)\big)$.
The graph $X$ is endowed with the graph distance. 
For a topological space $\T$, we denote by $C_\infty(\T)$
the space of all functions $f : \T \LA \mathbb{R}$ vanishing at infinity
with the uniform topology $\|f\|_\infty^{\T}$.

The ($p$-)transition operator $L$ acting on $C_\infty(X)$ is defined by
$$
Lf(x):=\sum_{e \in E_x}p(e)f\big(t(e)\big) \qquad (x \in V).
$$ 
The $n$-step transition probability $p(n, x, y) \, (n \in \mathbb{N}, \, x, y \in V)$
is given by $p(n, x, y):=L^n \delta_y(x)$, 
where $\delta_y( \cdot )$ stands for the Dirac delta function with pole at $y$. 
If there is a positive function $m : V \LA (0, \infty)$ up to constant multiple such that
$$
p(e)m\big(o(e)\big)=p(\ol{e})m\big(t(e)\big) \qquad (e \in E),
$$
the random walk is called ($m$-){\it symmetric} or {\it reversible}.  
Otherwise, it is called ($m$-){\it non-symmetric}. 

As one of the most fundamental examples of infinite graphs, 
a {\it crystal lattice} has been studied by many authors from both 
geometric and probabilistic viewpoints. 
Roughly speaking, an infinite graph $X$ is called a ($\Gamma$-)crystal lattice if  
$X$ is an infinite-fold covering graph of a finite graph 
whose covering transformation group $\Gamma$
is abelian. 
Typical examples we have in mind are
the square lattice, the triangular lattice and the hexagonal lattice, and so on.
For basic results, see \cite{Kotani, KS00, KS06, KSS} and literatures therein. 

In the theory of random walks on infinite graphs, to investigate
the long time asymptotics, for instance, 
the {\it central limit theorem} (CLT) is a principal topic 
for both geometers and probabilists. 
Recently, Ishiwata, Kawabi and Kotani \cite{IKK} proved 
two kinds of functional CLTs for non-symmetric random walks 
on crystal lattices using the theory of {\it discrete geometric analysis}
developed by Kotani and Sunada. 
For more details on discrete geometric analysis, see Section 2.
We also refer to \cite{Sunada, KS06}.

Before stating our results, 
we start with a brief review of the setting and results in \cite{IKK}. 
Let us consider a ($\Gamma$-)crystal lattice $X=(V, E)$, where 
the covering transformation group $\Gamma$, acting on $X$ freely, is a torsion free, finitely generated 
abelian group. Here we may assume that $\Gamma$ is isomorphic to $\mathbb{Z}^d$, 
without loss of generality. 
We denote by $X_0=(V_0, E_0)$ its (finite) quotient graph $\Gamma \backslash X$. 
Let $p : E \LA (0, 1]$ be a $\Gamma$-invariant transition probability. Namely, 
it satisfies (\ref{trans-prob})
and $p(\sigma e)=p(e)$ for every $\sigma \in \Gamma, \, e \in E$. 
Through the covering map $\pi : X \LA X_0$, the transition probability $p$ also
induces a Markov chain with values in $X_0$.  
Since the transition probability $p$ on $X$ is positive, 
the random walk on $X$ is {\it irreducible}, 
that is, for every $x, y \in V$, there exists $n=n(x, y) \in \mathbb{N}$ such that
$p(n, x, y)>0$. 
And so is the random walk on $X_0$.
Then, thanks to the Perron--Frobenius theorem, we find a unique {\it invariant probability
measure} on $V_0$ with
\begin{equation}\label{invariant}
\sum_{x \in V_0}m(x)=1 \quad \text{and} \quad 
m(x)=\sum_{e \in (E_0)_x}p(\ol{e})m\big(t(e)\big) \qquad (x \in V_0).
\end{equation}
It is also called a {\it stationary distribution} (see e.g., Durrett \cite{D}).
We also write $m : V \LA (0, 1]$ for the $\Gamma$-invariant lift of $m$ to $X$. 

Let $\h_1(X_0, \mathbb{R})$ and $\h^1(X_0, \mathbb{R})$ be
the first homology group and the first cohomology group of $X_0$, respectively.
We take a linear map $\rho_{\mathbb{R}}$ 
from $\h_1(X_0, \mathbb{R})$ onto $\Gamma \otimes \mathbb{R}$
through the covering map $\pi : X \LA X_0$. We define the {\it homological direction}
of the random walk on $X_0$ by
$$
\gamma_p:=\sum_{e \in E_0}p(e)m\big(o(e)\big)e \in \h_1(X_0, \mathbb{R}),
$$
and we call $\rho_{\mathbb{R}}(\gamma_p)(\in \Gamma \otimes \mathbb{R})$ 
the {\it asymptotic direction}. 
We remark that the random walk is $m$-symmetric if and only if $\gamma_p=0$.
Moreover, $\gamma_p=0$ implies $\rho_{\mathbb{R}}(\gamma_p)=\bm{0}$.
However, the converse does not hold in general. 
We write $g_0$ for the ($p$-){\it Albanese metric} on $\Gamma \otimes \mathbb{R}$.
(See Section 2, for its precise definition.) 
We call that a periodic realization 
$\Phi_0 : X \LA \Gamma \otimes \mathbb{R}$
is ($p$-){\it modified harmonic} if
\begin{equation} \label{modified harmonicity}
\sum_{e \in E_x}p(e) \Big( \Phi_0\big(t(e)\big) - \Phi_0\big(o(e)\big) \Big)
=\rho_{\mathbb{R}}(\gamma_p) \qquad (x \in V).
\end{equation}
This notion was first proposed in \cite{KS06} to seek the most canonical periodic 
realization of a topological crystal in the geometric context.  

Now we are in a position to review two kinds of CLTs formulated in \cite{IKK}.
We first set a reference point $x_* \in V$ with $\Phi_0(x_*)=\bm{0}$ and put
$\xi_n(c):=\Phi_0\big(w_n(c)\big) \, \big(n=0, 1, 2, \dots, \, c \in \Omega_{x_*}(X)\big)$.
Let $C_{\bm{0}}\big([0, \infty), (\Gamma \otimes \mathbb{R}, g_0)\big)$ be the set
of all continuous paths from $[0, \infty)$ to $(\Gamma \otimes \mathbb{R}, g_0)$ 
starting from the origin. 
We equip it with the usual compact uniform topology. 
We define a measurable map
$\X^{(n)} : \Omega_{x_*}(X) \LA 
\big( C_{\bm{0}}\big([0, \infty), (\Gamma \otimes \mathbb{R}, g_0)\big), \mu\big)$ by
$$
\X_t^{(n)}(c):=\frac{1}{\sqrt{n}}\Big\{ \xi_{[nt]}(c)+(nt-[nt])\big(\xi_{[nt]+1}(c)-\xi_{[nt]}(c)\big)
-nt \rho_{\mathbb{R}}(\gamma_p)\Big\} \qquad (t \geq 0),
$$
where $\mu$ is the Wiener measure on the path space
$C_{\bm{0}}\big([0, \infty), (\Gamma \otimes \mathbb{R}, g_0)\big)$.
We write $\Prob^{(n)} \, (n=1, 2, \dots)$ for the probability measure 
on $C_{\bm{0}}\big([0, \infty), (\Gamma \otimes \mathbb{R}, g_0)\big)$ induced by $\X^{(n)}$. 
Then the CLT of the first kind is stated as follows:

\begin{Thm} {\bf (\cite[Theorem 2.2]{IKK})}\label{CLT1}
The sequence $\{\Prob^{(n)}\}_{n=1}^\infty$ of probability measures 
converges weakly to the Wiener measure $\mu$ as $n \to \infty$. 
In other words, the sequence $\{\X^{(n)}\}_{n=1}^\infty$
converges to a $(\Gamma \otimes \mathbb{R}, g_0)$-valued standard 
Brownian motion $(B_t)_{t \geq 0}$ starting from the origin in law. 
\end{Thm}

Next we introduce a family $\{p_\ve\}_{0 \leq \ve \leq 1}$ of transition probabilities on $X$
by $p_\ve(e):=p_0(e)+\ve q(e) \, (e \in E)$, where 
$$
p_0(e):=\frac{1}{2}\Big( p(e)+\frac{m\big(t(e)\big)}{m\big(o(e)\big)}p(\ol{e})\Big), \quad 
q(e):=\frac{1}{2}\Big( p(e)-\frac{m\big(t(e)\big)}{m\big(o(e)\big)}p(\ol{e})\Big) \qquad (e \in E).
$$
This is nothing but the interpolation of the original transition probability $p=p_1$
and the $m$-symmetric transition probability $p_0$. 
We should note that, in this setting, the relation
 $\rho_{\mathbb{R}}(\gamma_{p_\ve})=\ve \rho_{\mathbb{R}}(\gamma_p)$ plays
 a crucial role to obtain the CLT of the second kind. 
We write $g_0^{(\ve)}$ for the ($p_\ve$-)Albanese metric 
on $\Gamma \otimes \mathbb{R}$. Moreover,
let $\Phi_0^{(\ve)} : X \LA (\Gamma \otimes \mathbb{R}, g_0^{(0)})$ be 
the $p_\ve$-modified harmonic realization of $X$. 

Now set a reference point $x_* \in V$ satisfying $\Phi_0^{(\ve)}(x_*)=\bm{0}$ 
for all $0 \leq \ve \leq 1$ and
put $\xi_n^{(\ve)}(c):=\Phi_0^{(\ve)}\big(w_n(c)\big) 
\, \big(n=0, 1, 2, \dots, c \in \Omega_{x_*}(X)\big)$. 
We define a measurable map $\Y^{(\ve, n)} : \Omega_{x_*}(X) 
\LA C_{\bm{0}}\big([0, \infty), (\Gamma \otimes \mathbb{R}, g_0^{(0)})\big)$ by
$$
\Y_t^{(\ve, n)}(c):=\frac{1}{\sqrt{n}}\Big\{\xi_{[nt]}^{(\ve)}(c)
+(nt-[nt])\big( \xi_{[nt]+1}^{(\ve)}(c) - \xi_{[nt]}^{(\ve)}(c)\big)\Big\} \qquad (t \geq 0).
$$
Let $\nu$ be the probability measure on 
$C_{\bm{0}}\big([0, \infty), (\Gamma \otimes \mathbb{R}, g_0^{(0)})\big)$
induced by the stochastic process 
$\big(B_t+\rho_{\mathbb{R}}(\gamma_p) t\big)_{t \geq 0}$, where
$(B_t)_{t \geq 0}$ is a $(\Gamma \otimes \mathbb{R}, g_0^{(0)})$-valued standard 
Brownian motion with $B_0=\bm{0}$.
 Moreover let $\Probq^{(\ve, n)}$ be the probability measure on 
 $C_{\bm{0}}\big([0, \infty), (\Gamma \otimes \mathbb{R}, g_0^{(0)})\big)$
 induced by $\Y^{(\ve, n)}$.
Then the CLT of the second kind is the following.

\begin{Thm} {\bf (\cite[Theorem 2.4]{IKK})} \label{CLT2}
The sequence $\{\Probq^{(n^{-1/2}, n)}\}_{n=1}^\infty$ of probability measures 
converges weakly to the probability measure $\nu$ as $n \to \infty$. 
In other words, the sequence $\{\Y^{(n^{-1/2}, n)}\}_{n=1}^\infty$ 
converges to a $(\Gamma \otimes \mathbb{R}, g_0^{(0)})$-valued standard 
Brownian motion with drift $\rho_{\mathbb{R}}(\gamma_p)$ starting from the origin in law. 
\end{Thm}

The main purpose of the present paper is to establish yet another CLT
for non-symmetric random walks on crystal lattices.
Our approach is inspired by a {\it measure-change techinique} due to 
Alexopoulos \cite{A} in which several limit theorems for random walks
on discrete groups of polynomial volume growth are obtained.

Now consider the ($m$-)non-symmetric transition probability $p : E \LA (0, 1]$.
In particular, we assume that $\rho_{\mathbb{R}}(\gamma_p) \neq \bm{0}$. 
Let $\Phi_0 : X \LA (\Gamma \otimes \mathbb{R}, g_0)$ be a 
($p$-)modified harmonic realization of $X$. 
We define a function 
$F=F_x(\lambda) : V_0 \times \Hom(\Gamma, \mathbb{R}) \LA (0, \infty)$ 
by
\begin{equation}\label{F}
F_x(\lambda):=\sum_{e \in (E_0)_x}p(e)
\exp\Big( {}_{\Hom(\Gamma, \mathbb{R})}
\big\la \lambda, \Phi_0\big(t(\widetilde{e})\big) - \Phi_0\big(o(\widetilde{e})\big)
\big\ra_{\Gamma \otimes \mathbb{R}}\Big),
\end{equation}
for $x \in V_0, \, \lambda \in \Hom(\Gamma, \mathbb{R})$,
where $\widetilde{e}$ stands for a lift of $e \in E_0$ to $X$.
Then we can verify that, for every $x \in V_0$, 
the function $F_x( \cdot) : \Hom(\Gamma, \mathbb{R}) \LA (0, \infty)$ has a unique 
minimizer $\lambda_*=\lambda_*(x) \in \Hom(\Gamma, \mathbb{R})$. 
(See Lemma \ref{key-lem}.) 
We define a positive function $\p : E_0 \LA (0, 1]$ by 
\begin{equation}\label{new-transition}
\p(e):=\frac{p(e)\exp\Big( {}_{\Hom(\Gamma, \mathbb{R})}
\big\la \lambda_*\big(o(e)\big), \Phi_0\big(t(\widetilde{e})\big) - \Phi_0\big(o(\widetilde{e})\big)
\big\ra_{\Gamma \otimes \mathbb{R}}\Big)}{F_{o(e)}\big(\lambda_*(o(e))\big)} \qquad (e \in E_0).
\end{equation}
Then it is straightforward to check that the function $\p$ 
also gives a transition probability on $X_0$. 
Noting the random walk $\{w_n^{(\p)}\}_{n=0}^\infty$ 
associated with $\p$ is also irreducible, 
the Perron--Frobenius theorem yields a unique normalized 
invariant measure $\m : V_0 \LA (0, 1]$ in the sense of (\ref{invariant}). 
We write $\p : E \LA (0, 1]$ and $\m : V \LA (0, 1]$ for the $\Gamma$-invariant
lifts of $\p : E_0 \LA (0, 1]$ and $\m : V_0 \LA (0, 1]$, respectively. 
We write $g_0^{(\p)}$ for the ($\p$-)Albanese metric associated with the
transition probability $\p$.

Let $L_{(\p)}$ be the transition operator, acting on $C_\infty(X)$, associated with 
the transition probability $\p$. Namely,
$$
L_{(\p)}f(x)=\sum_{e \in E_x}\p(e)f\big(t(e)\big) \qquad (x \in V).
$$
Recalling that the function $F=F_x(\lambda)$ has the (unique) minimizer 
$\lambda_*=\lambda_*(x)$ for every $x \in V_0$,
it follows that 
\begin{equation}
L_{(\p)}\Phi_0(x) - \Phi_0(x)=\sum_{e \in E_x}\p(e)
\Big( \Phi_0\big(t(e)\big)-\Phi_0\big(o(e)\big)\Big)=\bm{0} \qquad (x \in V).
\end{equation} 
This equation means that the $p$-modified harmonic realization 
$\Phi_0 : X \LA \Gamma \otimes \mathbb{R}$ is a $\p$-harmonic realization. 
In particular, we obtain $\rho_{\mathbb{R}}(\gamma_{\p})=\bm{0}$. 
Here we should emphasize that the transition probability $\p : E_0 \LA (0, 1]$ 
coincides with the original one $p : E_0 \LA (0, 1]$ 
provided that $\rho_{\mathbb{R}}(\gamma_p)=\bm{0}$. 

We fix a reference point $x_* \in V$ such that $\Phi_0(x_*)=\bm{0}$ and put 
$$
\xi_n^{(\p)}(c):=\Phi_0\big(w_n^{(\p)}(c)\big) \qquad
 \big(n=0, 1, 2, \dots, \, c \in \Omega_{x_*}(X)\big).
$$
We define a measurable map 
$\frak{X}^{(n)} : \Omega_{x_*}(X) 
\LA \big( C_{\bm{0}}\big([0, \infty), (\Gamma \otimes \mathbb{R}, g_0^{(\p)})\big), \mu \big)$ by
\begin{equation}\label{X}
\frak{X}_t^{(n)}(c):=\frac{1}{\sqrt{n}}\Big\{ \xi_{[nt]}^{(\p)}(c)+(nt-[nt])
\big(\xi_{[nt]+1}^{(\p)}(c)-\xi_{[nt]}^{(\p)}(c)\big)
\Big\} \qquad (t \geq 0),
\end{equation}
where $\mu=\mu^{(\p)}$ is the Wiener measure on 
$C_{\bm{0}}\big([0, \infty), (\Gamma \otimes \mathbb{R}, g_0^{(\p)})\big)$. 
We write $\frak{P}^{(n)} \, (n=1, 2, \dots)$ for the probability measure 
on $C_{\bm{0}}\big([0, \infty), (\Gamma \otimes \mathbb{R}, g_0^{(\p)})\big)$
 induced by $\frak{X}^{(n)}$. 
Then our main theorem is stated as follows:

\begin{Thm} \label{CLT3}
The sequence $\{\frak{P}^{(n)}\}_{n=1}^\infty$ of probability measures 
converges weakly to the Wiener measure $\mu$ as $n \to \infty$. 
Namely, the sequence $\{\frak{X}^{(n)}\}_{n=1}^\infty$  
converges to a $(\Gamma \otimes \mathbb{R}, g_0^{(\p)})$-valued standard 
Brownian motion $(B_t^{(\p)})_{t \geq 0}$ starting from the origin in law. 
\end{Thm}

Finally, we give a relationship between the $n$-step transition probabilities
$p(n, x, y)$ and $\p(n, x, y)$ as follows: 

\begin{Thm}\label{asymptotic-p}
There exist some positive constants $C_1$ and $C_2$ such that 
$$
C_1p(n, x, y)
\exp\big( n M_p\big) \leq \p(n, x, y) \leq C_2p(n, x, y)
\exp\big( n M_p\big)
$$
for all $n \in \mathbb{N}$ and $x, y \in V$, where
$$
M_p:=\sum_{\bm{x} \in V_0}m(\bm{x})\Big( 
{}_{\Hom(\Gamma, \mathbb{R})}\big\la\lambda_*(\bm{x}), 
\rho_{\mathbb{R}}(\gamma_p) \big\ra_{\Gamma \otimes \mathbb{R}}
- \log F_{\bm{x}}\big(\lambda_*(\bm{x})\big)\Big).
$$
\end{Thm}
For the precise long time asymptotic behavior of $p(n, x, y)$,
we refer to \cite[Theorem~2.5]{IKK} and Trojan \cite{Trojan}.

The rest of the present paper is organized as follows: 
In Section 2, we provide a brief review on the theory of 
discrete geometric analysis. In Section 3, we state our measure-change technique in details
and give proofs of Theorems \ref{CLT3} and \ref{asymptotic-p}.
We also discuss a relationship between our measure-change technique and 
a discrete analogue of Girsanov's theorem due to Fujita \cite{Fujita}. 
Finally, in Section 4, we give some concrete examples of non-symmetric random walks 
on crystal lattices. 

\section{A quick review on discrete geometric analysis}
In this section, we give basic materials of the theory of discrete geometric analysis
on graphs quickly. For more details, we refer to Kotani--Sunada \cite{KS06}
and Sunada \cite{Sunada}.

We consider a random walk on a finite graph $X_0=(V_0, E_0)$
associated with a transition probability $p : E_0 \LA (0, 1]$. 
Thanks to the Perron--Frobenius theorem, there is a unique (normalized) invariant measure 
$m : V_0 \LA (0, 1]$ in the sense of (\ref{invariant}). 
The random walk is called ($m$-){\it symmetric} if
$p(e)m\big(o(e)\big)=p(\ol{e})m\big(t(e)\big) \, (e \in E_0)$. 

First we define the 0-chain group and the 1-chain group of $X_0$ by
$$
C_0(X_0, \mathbb{R}):=\Big\{ \sum_{x \in V_0} a_x x \, \Big| \, a_x \in \mathbb{R} \Big\},
\quad
C_1(X_0, \mathbb{R}):=\Big\{ \sum_{e \in E_0} a_e e \, \Big| \, a_e \in \mathbb{R}
, \, \ol{e}=-e \Big\},
$$
respectively.
Let $\del : C_1(X_0, \mathbb{R}) \LA C_0(X_0, \mathbb{R})$ be the boundary map,
given by the homomorphism satisfying $\del(e):=t(e)-o(e) \, (e \in E_0)$.
The first homology group $\h_1(X_0, \mathbb{R})$ of $X_0$ is defined by
$\Ker(\del)\big( \subset C_1(X_0, \mathbb{R})\big)$.

On the other hand, we define the 0-cochain group and the 1-cochain group of $X_0$ by
$$
C^0(X_0, \mathbb{R}):=\{f : V_0 \LA \mathbb{R}\}, \quad 
C^1(X_0, \mathbb{R}):=\{\omega : E_0 \LA \mathbb{R} \, | \, \omega(\ol{e})=-\omega(e)\},
$$
respectively. The difference operator $d : C^0(X_0, \mathbb{R}) \LA C^1(X_0, \mathbb{R})$
is defined by the homomorphism with $df(e):=f\big(t(e)\big)-f\big(o(e)\big) \, (e \in E_0)$.
We also define $\h^1(X_0, \mathbb{R}):=C^1(X_0, \mathbb{R})/\mathrm{Im}\,(d)$, 
called the first cohomology group of $X_0$.

Next we define the transition operator 
$L : C^0(X_0, \mathbb{R}) \LA C^0(X_0, \mathbb{R})$ by
$$
Lf(x):=(I-\delta_pd)f(x)=
\sum_{e \in (E_0)_x}p(e)f\big(t(e)\big) \qquad \big(x \in V_0, \, f \in C^1(X_0, \mathbb{R})\big), 
$$
where the operator $\delta_p : C^1(X_0, \mathbb{R}) \LA C^0(X_0, \mathbb{R})$ is defined by
$$
\delta_p\omega(x):=-\sum_{e \in (E_0)_x}p(e)\omega(e)
 \qquad \big(x \in V_0, \, \omega \in C^1(X_0, \mathbb{R})\big).
$$
We introduce the quantity $\gamma_p$, called the {\it homological direction} 
of the given random walk on $X_0$, by 
$$
\gamma_p:=\sum_{e \in E_0}\widetilde{m}(e)e \in C_1(X_0, \mathbb{R}),
$$
where $\widetilde{m}(e):=p(e)m\big(o(e)\big) \, (e \in E_0)$. 
It is easy to show that $\del(\gamma_p)=0$, that is, $\gamma_p \in \h_1(X_0, \mathbb{R})$.
It should be noted that the transition probability $p$ gives an $m$-symmetric 
random walk on $X_0$ if and only if $\gamma_p=0$. 
A 1-form $\omega \in C^1(X_0, \mathbb{R})$ is said to be {\it modified harmonic} if
$$
\delta_p\omega(x)+
{}_{C_1(X_0, \mathbb{R})}\la \gamma_p, \omega \ra_{C^1(X_0, \mathbb{R})}
=0 \qquad (x \in V_0).
$$
We denote by $\mathcal{H}^1(X_0)$ 
the space of modified harmonic 1-forms on $X_0$, and equip it 
with the inner product
$$
\La \omega, \eta \Ra_p:=\sum_{e \in E_0}\widetilde{m}(e)\omega(e)\eta(e)
-\la \gamma_p, \omega \ra \la \gamma_p, \eta \ra 
\qquad \big(\omega, \eta \in \mathcal{H}^1(X_0)\big).
$$
Due to the discrete analogue of Hodge--Kodaira theorem (cf.~\cite[Lemma 5.2]{KS06}), 
we may identify 
$\big(\mathcal{H}^1(X_0), \La \cdot, \cdot \Ra_p\big)$ with $\h^1(X_0, \mathbb{R})$.

Now let $X=(V, E)$ be a $\Gamma$-crystal lattice. 
Namely, $X$ is a covering graph of a finite graph $X_0$
with an abelian covering transformation group $\Gamma$. 
We write $p : E \LA (0, 1]$ and $m : V \LA (0, 1]$ 
for the $\Gamma$-invariant lifts of $p : E_0 \LA (0, 1]$ and $m : V_0 \LA (0, 1]$, respectively. 
Through the covering map $\pi : X \LA X_0$, we take the surjective linear map
$\rho_{\mathbb{R}} : \h_1(X_0, \mathbb{R}) 
\LA \Gamma \otimes \mathbb{R}(\cong \mathbb{R}^d)$.
We consider the transpose 
${}^t \rho_{\mathbb{R}} : \Hom(\Gamma, \mathbb{R}) \LA \h^1(X_0, \mathbb{R})$,
which is a injective linear map.
Here $\Hom(\Gamma, \mathbb{R})$ 
denotes the space of homomorphisms from $\Gamma$ into $\mathbb{R}$.
We induce a flat metric $g_0$ on the Euclidean space $\Gamma \otimes \mathbb{R}$ 
through the following diagram:

$$
\xymatrix{ 
 (\Gamma \otimes \mathbb{R}, g_0) \ar @{<<-}[r]^{\rho_{\mathbb{R}}} 
 \ar @{<->}[d]^{\mathrm{dual}} & \h_1(X_0, \mathbb{R}) 
 \ar @{<->}[d]^{\mathrm{dual}} &\\
 \Hom(\Gamma, \mathbb{R})  \ar @{^{(}->}[r]_{{}^t \rho_{\mathbb{R}}} 
 & \h^1(X_0, \mathbb{R})  &\hspace{-1.5cm} \cong 
 & \hspace{-1.5cm}\big(\mathcal{H}^1(X_0), \La \cdot , \cdot \Ra_p\big).}
$$
This metric $g_0$ is called the {\it Albanese metric} on $\Gamma \otimes \mathbb{R}$.

From now on, we realize the crystal lattice $X$ into the continuous model
$(\Gamma \otimes \mathbb{R}, g_0)$ in the following manner.
A {\it periodic realization} of $X$ into $\Gamma \otimes \mathbb{R}$ is defined 
by a piecewise linear map $\Phi : X \LA \Gamma \otimes \mathbb{R}$
with $\Phi(\sigma x)=\Phi(x)+\sigma \otimes 1 \, (\sigma \in \Gamma, \, x \in V)$. 
We introduce a special periodic realization 
$\Phi_0 : X \LA \Gamma \otimes \mathbb{R}$ by
\begin{equation}\label{Alb}
{}_{\Hom(\Gamma, \mathbb{R})}\big\la \omega, 
\Phi_0(x) \big\ra_{\Gamma \otimes \mathbb{R}}
=\int_{x_*}^x \widetilde{\omega} 
\qquad \big(x \in V, \, \lambda \in \Hom(\Gamma, \mathbb{R})\big),
\end{equation}
where $x_*$ is a fixed reference point satisfying $\Phi_0(x_*)=\bm{0}$ and
$\widetilde{\omega}$ is the lift of $\omega$ to $X$. 
Here
$$
\int_{x_*}^x \widetilde{\omega}=\int_c \widetilde{\omega}
:=\sum_{i=1}^n \widetilde{\omega}(e)
$$
for a path $c=(e_1, \dots, e_n)$ with $o(e_1)=x_*$ and $t(e_n)=x$.
It should be noted that this line integral does not depend on the choice of a path $c$.
The periodic realization $\Phi_0$ given by above 
enjoys the so-called {\it modified harmonicity}
in the sense that
$$
L\Phi_0(x) - \Phi_0(x)=\rho_{\mathbb{R}}(\gamma_p) \qquad (x \in V).
$$
We note that this equation is also written as (\ref{modified harmonicity}).
Further, such a realization is uniquely determined up to translation. 
We call the quantity $\rho_{\mathbb{R}}(\gamma_p)$ the {\it asymptotic direction}
of the given random walk. 
We should emphasize that $\gamma_p=0$ implies
$\rho_{\mathbb{R}}(\gamma_p)=\bm{0}.$ 
However, the converse does not always hold, that is, there is a case 
$\gamma_p \neq 0$ and $\rho_{\mathbb{R}}(\gamma_p)=\bm{0}$.  
(See Subsection 4.2, for an example.)
If we equip $\Gamma \otimes \mathbb{R}$ with the Albanese metric, 
then the modified harmonic realization 
$\Phi _0 : X \LA (\Gamma \otimes \mathbb{R}, g_0)$
is said to be a {\it modified standard realization}.


\section{Proofs of the main results}

\subsection{A measure--change technique}

In what follows, we write $\lambda[\x]_{\Gamma \otimes \mathbb{R}}
:={}_{\Hom(\Gamma, \mathbb{R})}\la \lambda, \x \ra_{\Gamma \otimes \mathbb{R}}$
$(\lambda \in \Hom(\Gamma, \mathbb{R}), \, \x \in \Gamma \otimes \mathbb{R})$ and
$$
d\Phi_0(e):=\Phi_0\big(t(e)\big) - \Phi_0\big(o(e)\big) \qquad (e \in E),
$$
for brevity. Take an orthonormal basis $\{\omega_1, \dots, \omega_d\}$
in $\Hom(\Gamma, \mathbb{R})
\big(\subset (\mathcal{H}^1(X_0), \La \cdot, \cdot \Ra_p)\big)$, 
and denote by 
$\{\vv_1, \dots, \vv_d\}$ its dual basis in $\Gamma \otimes \mathbb{R}$.
Namely, $\omega_i[\vv_j]_{\Gamma \otimes \mathbb{R}}=\delta_{ij} \, (1 \leq i, j \leq d)$. 
Then, $\{\vv_1, \dots, \vv_d\}$ is an orthonormal basis
in $\Gamma \otimes \mathbb{R}$ with respect to the Albanese metric $g_0$.
We may identify
$\lambda=\lambda_1\omega_1+ \dots +\lambda_d\omega_d 
\in \Hom(\Gamma, \mathbb{R})$ with
$(\lambda_1, \dots, \lambda_d) \in \mathbb{R}^d$.
Furthermore, we write $x_i:=\omega_i[\x]
_{\Gamma \otimes \mathbb{R}}$, 
$\Phi_0(x)_i:=\omega_i[\Phi_0(x)]_{\Gamma \otimes \mathbb{R}}$
and $\del_i:=\del/\del \lambda_i \, (i=1, \dots, d, \, x \in V)$.
We denote by $O(\cdot)$ the Landau symbol. 

At the beginning, consider the function 
$F=F_x(\lambda) : V_0 \times \Hom(\Gamma, \mathbb{R}) \LA \mathbb{R}$
defined by (\ref{F}).
We easily see that 
$F=F_x(\lambda)$ is a positive function on 
$V_0 \times \Hom(\Gamma, \mathbb{R})$ with $F_x(\bm{0})=1 \, (x \in V_0)$.  
In our setting, the following lemma plays a siginificant role 
so as to obtain Theorem \ref{CLT3}. 

\begin{Lem}\label{key-lem}
For every $x \in V_0$, the function 
$F_x( \cdot) : \Hom(\Gamma, \mathbb{R}) \LA (0, \infty)$
has a unique minimizer $\lambda_*=\lambda_*(x)$. 
\end{Lem}

\noindent
{\bf Proof.} 
Fix a fixed $x \in V_0$, we have 
$$
\begin{aligned}
\del_iF_x(\lambda)&=
\del_i \Bigg( \sum_{e \in (E_0)_x}p(e) 
\exp \Big( \lambda\big[ d\Phi_0(\widetilde{e})\big]
_{\Gamma \otimes \mathbb{R}}\Big) \Bigg)\\
&=\del_i \Bigg( \sum_{e \in (E_0)_x}p(e) 
\exp \Big( \sum_{i=1}^d \lambda_i \cdot \omega_i\big[ d\Phi_0(\widetilde{e})\big]
_{\Gamma \otimes \mathbb{R}}\Big) \Bigg)\\
&=\sum_{e \in (E_0)_x}p(e) 
\exp\Big( \lambda\big[ d\Phi_0(\widetilde{e})\big]_{\Gamma \otimes \mathbb{R}}\Big) 
d\Phi_0(\widetilde{e})_i  \qquad \big(i=1, \dots, d, \, \lambda \in \Hom(\Gamma, \mathbb{R})\big).
\end{aligned}
$$
In other words,
\begin{align}\label{nabla}
&\quad \Big( \del_1 F_x(\lambda), \dots, \del_d F_x(\lambda) \Big)\nn\\
&=\sum_{e \in (E_0)_x}p(e)
 \exp\Big( \lambda\big[ d\Phi_0(\widetilde{e})\big]
 _{\Gamma \otimes \mathbb{R}} \Big) 
 d\Phi_0(\widetilde{e}) \,(\in \Gamma \otimes \mathbb{R})
  \qquad \big(\lambda \in \Hom(\Gamma, \mathbb{R})\big).
\end{align}
Repeating the above calculation, we have 
$$
\del_i \del_j F_x(\lambda)=
\sum_{e \in (E_0)_x}p(e)  \exp\Big( \lambda\big[ d\Phi_0(\widetilde{e})\big]
 _{\Gamma \otimes \mathbb{R}} \Big)d\Phi_0(\widetilde{e})_i d\Phi_0(\widetilde{e})_j
$$ 
for $\lambda \in \Hom(\Gamma, \mathbb{R})$ and $1 \leq i, j \leq d$. 
Then, we know that $\big( \del_i \del_j F_x(\cdot) \big)_{1 \leq i, j \leq d} $, 
the {\it Hessian matrix} of
the function $F_x(\cdot)$, is positive definite. 
Indeed, consider the quadratic form corresponding to the Hessian matrix. Since 
\begin{align}\label{quadratic}
&\quad \sum_{1 \leq i, j \leq d}\sum_{e \in (E_0)_x}p(e) 
 \exp\Big( \lambda\big[ d\Phi_0(\widetilde{e})\big]
 _{\Gamma \otimes \mathbb{R}} \Big) 
 d\Phi_0(\widetilde{e})_i d\Phi_0(\widetilde{e})_j \xi_i \xi_j \nn\\
&= \sum_{e \in (E_0)_x}p(e) 
\exp\Big( \lambda\big[ d\Phi_0(\widetilde{e})\big]
 _{\Gamma \otimes \mathbb{R}}  \Big) 
 \Big\{ \sum_{i=1}^d d\Phi_0(\widetilde{e})_i  \xi_i \Big\}^2 
 \geq 0 \quad \big(\bm{\xi}=(\xi_1, \dots, \xi_d) \in \mathbb{R}^d\big)
\end{align}
and the transition probability $p$ is positive,
we easily see that the Hessian matrix is non-negative definite. 
By multiplying both sides of (\ref{quadratic}) by $m(x)$ and  
taking the sum which runs over all vertices of $X_0$, it readily follows that 
$$
\sum_{e \in E_0}\widetilde{m}(e)\exp\Big( \lambda\big[ d\Phi_0(\widetilde{e})\big]
 _{\Gamma \otimes \mathbb{R}}  \Big)
  \Big\{ \sum_{i=1}^d d\Phi_0(\widetilde{e})_i  \xi_i \Big\}^2 
  \geq 0, \qquad \big(\bm{\xi}=(\xi_1, \dots, \xi_d) \in \mathbb{R}^d\big).
$$
Next suppose that the left-hand side of (\ref{quadratic}) is zero. Then we have
$$
\sum_{i=1}^d d\Phi_0(\widetilde{e})_i \xi_i=0
$$
for all $e \in E_0$. 
This equation implies
 $\la \Phi_0(x), \bm{\xi} \ra_{\mathbb{R}^d}=\la \Phi_0(y), 
 \bm{\xi} \ra_{\mathbb{R}^d}$ for all $x, y \in V$, where
$\la \cdot , \cdot \ra_{\mathbb{R}^d}$ stands for 
the standard inner product on $\mathbb{R}^d$. 
 Let $\sigma_1, \dots, \sigma_d$ be generators of $\Gamma \cong \mathbb{Z}^d$. 
It follows from the periodicity of $\Phi_0$ that
$\la \sigma_i, \bm{\xi} \ra_{\mathbb{R}^d}=0 \, (1 \leq i \leq  d)$. 
Hence we conclude $\bm{\xi}=\bm{0}$. Namely,  
we have proved the positive definiteness of the Hessian matrix. 

This implies that the function 
$F_x( \cdot ) : \Hom(\Gamma, \mathbb{R}) \LA (0, \infty)$ 
is strictly convex for every $x \in V_0$.
Moreover, it is easily observed that 
$$
\lim_{|\lambda|_{\mathbb{R}^d}\to \infty} F_x(\lambda)=\infty \qquad (x \in X_0),
$$
due to its definition. 
Consequently, we know that there exists a unique minimizer 
$\lambda_*=\lambda_*(x) \in \Hom(\Gamma, \mathbb{R})$
of $F_x(\lambda)$ for each $x \in V_0$, 
thereby completing the proof. \qed 

\vspace{2mm}

Now consider the positive function $\p : E_0 \LA (0, 1]$ given by 
(\ref{new-transition}). 
By definition, we easily see that 
the function $\p$ also gives a positive transition probability on $X_0$.
Thus the transition probability $\p : E_0 \LA (0, 1]$ yields an irreducible random walk
$(\Omega_{x_*}(X), \widehat{\mathbb{P}}_{x_*}, \{w_n^{(\p)}\}_{n=0}^\infty)$ with values in $X$. 
Applying the Perron-Frobenius theorem again, 
we find a unique positive function $\frak{m} : V_0 \LA (0, 1]$ satisfying (\ref{invariant}).
Put $\widetilde{\frak{m}}(e):=\p(e)\m(o(e)) \, (e \in E_0)$. 
We also denote by $\p : E \LA (0, 1]$ and 
$\m : V \LA (0, 1]$ the $\Gamma$-invariant lifts of $\p : E_0 \LA (0, 1]$ 
and $\m : V_0 \LA (0, 1]$, respectively. 
As in the previous section, 
we construct the ($\p$-)Albanese metric $g_0^{(\p)}$ on 
$\Gamma \otimes \mathbb{R}$ 
associated with the transition probability $\p$. 
We take an orthonormal basis $\{\omega_1^{(\p)}, \dots, \omega_d^{(\p)}\}$
in $\Hom(\Gamma, \mathbb{R})
\big(\subset (\mathcal{H}^1(X_0), \La \cdot, \cdot \Ra_\p)\big)$.

We introduce the transition operator $L_{(\p)} : C_\infty(X) \LA C_\infty(X)$ 
associated with the transition probability $\p$ by
$$
L_{(\p)}f(x):=\sum_{e \in E_x}\p(e)f(t(e)) \qquad (x \in V).
$$
Recalling (\ref{nabla}) and the definition of $\lambda_*=\lambda_*(x)$, 
we see that
\begin{align}
\quad \Big( \del_1 F_x\big(\lambda_*(x)\big), \dots, 
\del_d F_x\big(\lambda_*(x)\big) \Big)=\sum_{e \in (E_0)_x}p(e)
 \exp\Big( \lambda_*(x)\big[ d\Phi_0(\widetilde{e})\big]
 _{\Gamma \otimes \mathbb{R}} \Big) 
 d\Phi_0(\widetilde{e})=\bm{0}\nn
\end{align}
holds for every $x \in V_0$. 
This immediately leads to
\begin{equation}\label{harmonicity}
L_{(\p)}\Phi_0(x)-\Phi_0(x)=\sum_{e \in E_x}\p(e)d\Phi_0(e)=\bm{0} \qquad (x \in V).
\end{equation}
From this equation, one concludes that the given $p$-modified standard realization
$\Phi_0 : X \LA (\Gamma \otimes \mathbb{R}, g_0)$
in the sense of (\ref{modified harmonicity})
is the harmonic realization associated with the changed transition probability $\p$.

\begin{Rem}
Equation 
{\rm (\ref{harmonicity})} implies $\rho_{\mathbb{R}}(\gamma_\p)=\bm{0}$. 
We also emphasize that the transition probability $\p : E_0 \LA (0, 1]$ 
coincides with the original one $p : E_0 \LA (0, 1]$ 
provided that $\rho_{\mathbb{R}}(\gamma_p)=\bm{0}$.  
\end{Rem}

\begin{Rem}
In our setting, it is essential to assume 
that the given transition probability $p$ is positive. 
Because, if it were not for the positivity of $p$, 
the assertion of {\rm Lemma \ref{key-lem}} would not hold in general. 
(There is a case where the function $F_x(\cdot)$ has no minimizers.)
On the other hand, to obtain {\rm Theorems \ref{CLT1}} and {\rm \ref{CLT2}}, 
it is sufficient to impose that
the given transition probability $p$ is non-negative with 
$p(e)+p(\ol{e})>0 \, (e \in E)$.  
\end{Rem}

\subsection{Proofs of Theorems \ref{CLT3} and \ref{asymptotic-p}}

This subsection is devoted to proofs of Theorems \ref{CLT3} and \ref{asymptotic-p}. 
Following the argument as in \cite[Theorem 2.2]{IKK} for the 
random walk associated with the changed transition probability $\p$,
we can carry out the proof of Theorem \ref{CLT3}.
Though a minor change of the proof is required, the argument is a little bit easier
due to $\rho_{\mathbb{R}}(\gamma_\p)=\bm{0}$.

As the first step, we prove the following lemma.

\begin{Lem}\label{proof-lem1}
For any $f \in C_0^\infty(\Gamma \otimes \mathbb{R})$,
as $N \to \infty$, $\ve \searrow 0$ and $N^2\ve \searrow 0$, we have
$$
\Big\|\frac{1}{N\ve^2}(I-L_{(\p)}^N)P_\ve f - P_\ve \Big(\frac{\Delta_{(\p)}}{2}f\Big)\Big\|_\infty^X \LA 0.
$$
Here $P_\ve : C_\infty(\Gamma \otimes \mathbb{R}) \LA C_\infty(X) \, (0 \leq \ve \leq 1)$
 is a scaling operator defined by
$$
P_\ve f(x):=f\big( \ve \Phi_0(x)\big) \qquad (x \in X)
$$
and $\Delta_{(\p)}$ stands for the positive Laplacian $-\sum_{i=1}^d(\del^2/\del x_i^2)$ on $\Gamma \otimes \mathbb{R}$
associated with the $\p$-Albanese metric $g_0^{(\p)}$.
\end{Lem}

\noindent
{\bf Proof.} We define the function 
$A^N(\Phi_0)_{ij} : V \LA \mathbb{R} \, (i, j=1, \dots, d, \, N \in \mathbb{N})$
by 
$$
A^N(\Phi_0)_{ij}(x):=\sum_{c \in \Omega_{x, N}(X)}
\p(c)\Big(\Phi_0\big(t(c)\big) - \Phi_0(x)\Big)_i
\Big(\Phi_0\big(t(c)\big) - \Phi_0(x)\Big)_j \qquad (x \in V),
$$
where $\p(c):=\p(e_1) \cdots \p(e_N)$ for $c=(e_1, \dots, e_N) \in \Omega_{x, N}(X)$. 
Applying Taylor's expansion formula, we have
\begin{align}\label{Taylor}
(I-L_{(\p)}^N)P_\ve f(x)&=-\ve \sum_{i=1}^d \frac{\del f}{\del x_i}\big( \ve \Phi_0(x)\big)
\sum_{c \in \Omega_{x, N}(X)}\p(c)\Big(\Phi_0\big(t(c)\big) - \Phi_0(x)\Big)_i \nn \\
&\qquad -\frac{\ve^2}{2}\sum_{1 \leq i, j \leq d}\frac{\del^2f}{\del x_i\del x_j}\big( \ve \Phi_0(x)\big)
A^N(\Phi_0)_{ij}(x)+O\big((N\ve)^3\big).
\end{align}
We see that the first term of the right-hand side of (\ref{Taylor}) equals 0 
due to the $\p$-harmonicity of $\Phi_0$. Next we define the function 
$\A(\Phi_0)_{ij} : V_0 \LA \mathbb{R} \, (i, j=1, \dots, d)$
by 
$$
\A(\Phi_0)_{ij}(x):=\sum_{e \in (E_0)_x}
\p(e)d\Phi_0(\widetilde{e})_i
d\Phi_0(\widetilde{e})_j \qquad (x \in V).
$$
We note that 
$\A(\Phi_0)_{ij}\big(\pi(x)\big)=A^1(\Phi_0)_{ij}(x) \, (x \in V, \, i, j=1, \dots, d)$
because $A^N(\Phi_0)_{ij}$ is $\Gamma$-invariant.
Then, using the $\p$-harmonicity again, 
$$
A^N(\Phi_0)_{ij}(x)=\sum_{k=0}^{N-1}
L_{(\p)}^k \big( \A(\Phi_0)_{ij}\big)\big(\pi(x)\big) \qquad (x \in V).
$$
Applying the ergodic theorem for $L_{(\p)}$ (cf. \cite[Theorem 3.2]{IKK}), we have
$$
\frac{1}{N}\sum_{k=0}^{N-1}
L_{(\p)}^k \big( \A(\Phi_0)_{ij}\big)\big(\pi(x)\big)
=\sum_{x \in V_0}\m(x) \A(\Phi_0)_{ij}(x) + O\Big(\frac{1}{N}\Big).
$$
Then, (\ref{Alb}) and $\p$-harmonicity of $\Phi_0$ imply
$$
\sum_{x \in V_0}\m(x) \A(\Phi_0)_{ij}(x)=\sum_{e \in E_0}\widetilde{\m}(e)
\omega_i^{(\p)}(e)\omega_j^{(\p)}(e)=\La \omega_i^{(\p)}, \omega_j^{(\p)}\Ra_{\p}
=\delta_{ij}
$$
for $1 \leq i, j \leq d$. Putting it all together, we obtain
$$
\frac{1}{N\ve^2}(I-L_{(\p)}^N)P_\ve f =
P_\ve \Big(\frac{\Delta_{(\p)}}{2}f\Big)+O(N^2\ve)+O\Big(\frac{1}{N}\Big).
$$
Finally, letting $N \to \infty$, $\ve \searrow 0$ and $N^2\ve \searrow 0$, 
we complete the proof. \qed

\vspace{2mm}
Lemma \ref{proof-lem1} immediately leads to the following lemma. 
(See \cite[Theorem 2.1 and Lemma 4.2]{IKK} for details.)

\begin{Lem}\label{proof-lem2}
{\rm (1)} 
For any $f \in C_\infty(\Gamma \otimes \mathbb{R})$,
and $0 \leq s \leq t$, we have
$$
\lim_{n \to \infty}\Big\|L_{(\p)}^{[nt]-[ns]}P_{n^{-1/2}}f - P_{n^{-1/2}}\e^{-(t-s)\Delta_{(\p)}/2}f\Big\|_\infty^X=0.
$$

\noindent
{\rm (2)}
We fix $0 \leq t_1<\dots<t_\ell<\infty \, (\ell \in \mathbb{N})$.
Then, 
$$
(\frak{X}_{t_1}^{(n)}, \dots, \frak{X}_{t_\ell}^{(n)}) 
\overset{(d)}{\longrightarrow}
(B_{t_1}^{(\p)}, \dots, B_{t_\ell}^{(\p)}) \qquad (n \to \infty),
$$
where $(B_t^{(\p)})_{t \geq 0}$ is a 
$(\Gamma \otimes \mathbb{R}, g_0^{(\p)})$-valued standard 
Brownian motion with $B_{0}^{(\p)}=\bm{0}$.
\end{Lem}

Having obtained Lemma $\ref{proof-lem2}$,
it is sufficient to show the tightness of $\{\frak{P}^{(n)}\}_{n=1}^\infty$
for completing the proof of Theorem \ref{CLT3}.

\begin{Lem}\label{proof-lem4}
The sequence $\{\frak{P}^{(n)}\}_{n=1}^\infty$ is tight in 
$C_{\bm{0}}\big([0, \infty), (\Gamma \otimes \mathbb{R}, g_0^{(\p)})\big)$.
\end{Lem}

\noindent
{\bf Proof.} 
Throughout the proof, $C$ denotes a positive constant 
that may change at every occurrence. 
We put 
$\|d\Phi_0\|_\infty:=\max_{e \in E_0}\|d\Phi_0(\widetilde{e})\|_{g_0^{(\p)}}$. 

By virtue of the celebrated Kolmogorov's criterion, it is sufficient to show that
there exists some positive constant $C$ independent of $n$ such that
\begin{equation}\label{eq4}
\mathbb{E}^{\widehat{\mathbb{P}}_{x_*}}
\Big[ \big\| \frak{X}_t^{(n)}-\frak{X}_s^{(n)}\big\|_{g_0^{(\p)}}^4\Big]
\leq C(t-s)^2 \qquad (0 \leq s \leq t, \, n \in \mathbb{N}).
\end{equation}
We split the proof into two cases: {\bf (I)} : $t-s<n^{-1}$, {\bf (II)} : $t-s \geq n^{-1}$. 

First we consider the case {\bf (I)}. 
In both cases $ns \geq [nt]$ and $ns < [nt]$, we have
$$
\begin{aligned}
\big\| \frak{X}_t^{(n)}-\frak{X}_s^{(n)}\big\|_{g_0^{(\p)}}
&\leq n^{1/2}(t-s)\Big\{\big\|\xi_{[nt]+1}^{(\p)}-\xi_{[nt]}^{(\p)}\|_{g_0^{(\p)}}
+\big\|\xi_{[nt]}^{(\p)}-\xi_{[nt]-1}^{(\p)}\|_{g_0^{(\p)}}\Big\}\\
&\leq 2\|d\Phi_0\|_\infty n^{1/2}(t-s).
\end{aligned}
$$
Noting $n^2(t-s)^2<1$, we obatin
the desired estimate (\ref{eq4}) for case {\bf (I)}. 

Next we consider the case {\bf (II)}. 
Let $\mathcal{F}$ be the fundamental domain in $X$ containing $x_* \in V$
and define $\frak{M}_i^\ell=\frak{M}_i^\ell(\Phi_0) : 
V \LA \mathbb{R} \, (i=1, \dots, d, \, \ell=1, 2, 3, 4)$ by
$$
\frak{M}_i^\ell(x):=\sum_{e \in E_x}\p(e)d\Phi_0(e)_i^\ell \qquad (x \in V).
$$
We note that $\frak{M}_i^\ell$ is $\Gamma$-invariant and 
$\|\frak{M}_i^\ell\|_\infty^X \leq \|d\Phi_0\|_\infty^\ell \, (i=1, \dots, d).$
Moreover, we obtain $\frak{M}_i^1 \equiv 0 \, (i=1, \dots, d)$ due to
the $\p$-harmonicity of $\Phi_0$. 

Here we give a bound on 
$\mathbb{E}^{\widehat{\mathbb{P}}_{x_*}}
\Big[ \big\| \mathcal{X}_{\frac{M}{n}}^{(n)}-\mathcal{X}_{\frac{N}{n}}^{(n)}
\big\|_{g_0^{(\p)}}^4\Big] \, (n \in \mathbb{N}, \, M \geq N \in \mathbb{N})$.
First of all, we have
\begin{equation}\label{eq5}
\mathbb{E}^{\widehat{\mathbb{P}}_{x_*}}
\Big[ \big\| \mathcal{X}_{\frac{M}{n}}^{(n)}-\mathcal{X}_{\frac{N}{n}}^{(n)}
\big\|_{g_0^{(\p)}}^4\Big]
\leq Cn^{-2}\max_{i=1, \dots, d}\max_{x \in \mathcal{F}}
\Big\{ \sum_{c \in \Omega_{x, M-N}(X)}
\p(c)\Big( \Phi_0\big(t(c)\big)-\Phi_0(x)\Big)_i^4\Big\}.
\end{equation}
Now fix $i=1, \dots, d$ and $x \in \mathcal{F}$.
For $k=1, \dots, M-N$, we have
\begin{align}\label{eq6}
&\sum_{c \in \Omega_{x, k}(X)}\p(c)\Big( \Phi_0\big(t(c)\big)-\Phi_0(x)\Big)_i^4 \nn\\
&=\sum_{c' \in \Omega_{x, k-1}(X)}\p(c')\sum_{e \in E_{t(c')}}\p(e)
\Big\{ \Big( \Phi_0\big(t(e)\big)-\Phi_0\big(o(e)\big)\Big)_i+
\Big( \Phi_0\big(o(e)\big)-\Phi_0(x)\Big)_i\Big\} \nn\\
&= \sum_{c' \in \Omega_{x, k-1}(X)}\p(c')\frak{M}_i^4\big(t(c')\big)
        +4\sum_{c' \in \Omega_{x, k-1}(X)}
        \p(c')\Big( \Phi_0\big(t(c')\big)-\Phi_0(x)\Big)_i\frak{M}_i^3\big(t(c')\big)\nn\\
        &\quad +6\sum_{c' \in \Omega_{x, k-1}(X)}
        \p(c')\Big( \Phi_0\big(t(c')\big)-\Phi_0(x)\Big)_i^2\frak{M}_i^2\big(t(c')\big)\nn\\
        &\quad 
        +\sum_{c' \in \Omega_{x, k-1}(X)}
        \p(c')\Big( \Phi_0\big(t(c')\big)-\Phi_0(x)\Big)_i^4\nn \\
&\leq \sum_{c' \in \Omega_{x, k-1}(X)}
        \p(c')\Big( \Phi_0\big(t(c')\big)-\Phi_0(x)\Big)_i^4  \nn \\
&\quad +\|\frak{M}_i^4\|_\infty^X
        +4\|d\Phi_0\|_\infty\|\frak{M}_i^3\|_\infty^X +6\|\frak{M}_i^2\|_\infty^X 
        \sum_{c \in \Omega_{x, k-1}(X)}
        \p(c)\Big( \Phi_0\big(t(c)\big)-\Phi_0(x)\Big)_i^2.
\end{align}
Moreover, the $\p$-harmonicity implies
\begin{align}\label{eq7}
&\sum_{c \in \Omega_{x, k-1}(X)}
        \p(c)\Big( \Phi_0\big(t(c)\big)-\Phi_0(x)\Big)_i^2\nn\\
&=\sum_{c' \in \Omega_{x, k-2}(X)}\p(c')\Big\{
        \Big( \Phi_0\big(t(c')\big)-\Phi_0(x)\Big)_i^2 + \frak{M}_i^2\big(t(c')\big)\Big\}
\leq (k-1)\|d\Phi_0\|_\infty^2.
\end{align}
It follows from (\ref{eq6}) and (\ref{eq7}) that
\begin{align}\label{eq8}
&\sum_{c \in \Omega_{x, k}(X)}\p(c)\Big( \Phi_0\big(t(c)\big)-\Phi_0(x)\Big)_i^4\nn \\
&\leq\sum_{c \in \Omega_{x, k-1}(X)}\p(c)\Big( \Phi_0\big(t(c)\big)-\Phi_0(x)\Big)_i^4 
+ 5\|d\Phi_0\|_\infty^4 + 6\|d\Phi_0\|_\infty^4(k-1)\leq Ck^2.
\end{align}
Putting $k=M-N$, $M=[nt]+i, \, N=[ns]+j \, (i, j=0, 1)$ and combining 
(\ref{eq5}) with (\ref{eq8}), we obtain
$$
\begin{aligned}
\mathbb{E}^{\widehat{\mathbb{P}}_{x_*}}
\Big[ \big\| \frak{X}_t^{(n)}-\frak{X}_s^{(n)}\big\|_{g_0^{(\p)}}^4\Big] &\leq 
\mathbb{E}^{\widehat{\mathbb{P}}_{x_*}}
\Big[ \max_{i, j=0, 1}\big\| \mathcal{X}_{\frac{[nt]+i}{n}}^{(n)}
-\mathcal{X}_{\frac{[ns]+j}{n}}^{(n)}\big\|_{g_0^{(\p)}}^4\Big] \\
&\leq Cn^{-2}\big( [nt]-[ns]+1\big)^2 \\
&\leq C\Big( t-s+\frac{2}{n}\Big)^2 \leq C\big\{(t-s)+2(t-s)\big\}^2 \leq C(t-s)^2,
\end{aligned}
$$
where we used $[nt]-[ns] \leq n(t-s)+1$ and $n^{-1} \leq t-s$. 
Therefore, we have shown the the desired estimate (\ref{eq4}) for case {\bf (II)}. \qed
\vspace{2mm}

Next we prove Theorem \ref{asymptotic-p}. 
\vspace{2mm}

\noindent
{\bf Proof of Theorem \ref{asymptotic-p}.} 
For $n \in \mathbb{N}$ and $x, y \in V$, we have
$$
\begin{aligned}
\p(n, x, y)&=\sum_{\substack{(e_1, \dots, e_n) 
\in \Omega_{x, n}(X) \\ o(e_1)=x, \, t(e_n)=y}}
\p(e_1) \cdots \p(e_n)\\
&=\sum_{\substack{(e_1, \dots, e_n) 
\in \Omega_{x, n}(X) \\ o(e_1)=x, \, t(e_n)=y}}
p(e_1) \cdots p(e_n) \cdot \exp\Big( \sum_{i=1}^n
\lambda_*\big(o(e_i)\big)[ d\Phi_0(\widetilde{e}_i)]_{\Gamma \otimes \mathbb{R}}
\Big)\\
& \hspace{2cm} \times F_{o(e_1)}\big( \lambda_*(o(e_1))\big)^{-1} \cdots 
F_{o(e_n)}\big( \lambda_*(o(e_n))\big)^{-1}\\
&=\sum_{\substack{(e_1, \dots, e_n) 
\in \Omega_{x, n}(X) \\ o(e_1)=x, \, t(e_n)=y}}
p(e_1) \cdots p(e_n) \\
& \hspace{2cm} \times \exp\Big( \sum_{i=1}^n
\lambda_*\big(o(e_i)\big)[ d\Phi_0(\widetilde{e}_i)]_{\Gamma \otimes \mathbb{R}}
 -\sum_{i=1}^n \log F_{o(e_i)}\big( \lambda_*(o(e_i))\big) \Big).
\end{aligned}
$$
Using the ergodic theorem for 1-chains (cf.\,\cite{KS06}):
\begin{equation}\label{ergodic}
\frac{1}{n}\sum_{i=1}^n f(e_i)
=\sum_{e \in E_0}\widetilde{m}(e)f(e) + O\Big(\frac{1}{n}\Big) \qquad (f : E_0 \LA \mathbb{R}),
\end{equation}
we obtain
$$
\begin{aligned}
&\quad\frac{1}{n}\sum_{i=1}^n 
\Big(\lambda_*\big(o(e_i)\big)[ d\Phi_0(\widetilde{e}_i)]_{\Gamma \otimes \mathbb{R}}
 -\log F_{o(e_i)}\big( \lambda_*(o(e_i))\big) \Big) \\
 &=\sum_{e \in E_0}\widetilde{m}(e)
 \Big(\lambda_*\big(o(e)\big)[ d\Phi_0(\widetilde{e})]_{\Gamma \otimes \mathbb{R}}
 -\log F_{o(e)}\big( \lambda_*(o(e))\big) \Big) +O\Big(\frac{1}{n}\Big)\\
 &=\sum_{\bm{x} \in V_0}m(\bm{x}) 
  \Big(\lambda_*(\bm{x})\Big[ \sum_{e \in (E_0)_{\bm{x}}}d\Phi_0(\widetilde{e})\Big]
  _{\Gamma \otimes \mathbb{R}}
 -\log F_{\bm{x}}\big( \lambda_*(\bm{x})\big) \Big) +O\Big(\frac{1}{n}\Big)\\
 &=\sum_{\bm{x} \in V_0}m(\bm{x}) 
  \Big(\lambda_*(\bm{x})\big[ \rho_{\mathbb{R}}(\gamma_p)\big]
  _{\Gamma \otimes \mathbb{R}}
 -\log F_{\bm{x}}\big( \lambda_*(\bm{x})\big) \Big) +O\Big(\frac{1}{n}\Big)
 \end{aligned}
$$
for $x, y \in V$. Here we used the $p$-modified harmonicity of $\Phi_0$ for the final line. 
Finally, we obtain
$$
\p(n, x, y) = p(n, x, y)\exp\Big( n\sum_{\bm{x} \in V_0}m(\bm{x})
\Big(\lambda_*(\bm{x})
\big[ \rho_{\mathbb{R}}(\gamma_p)\big]
_{\Gamma \otimes \mathbb{R}} - 
\log F_{\bm{x}}\big(\lambda_*(\bm{x})\big)\Big)+O(1)\Big)
$$ 
for $x, y \in V$. This completes the proof. \qed

\begin{Rem}
Let us consider a special case where the $\Gamma$-crystal lattice $X$ is given by
a covering graph of an $\ell$-bouquet graph $(\ell \in \mathbb{N})$ consisting of 
one vertex $\bm{x} \in V_0$ and $\ell$-loops. 
Without using the ergodic theorem {\rm (\ref{ergodic})} in the proof of {\rm Theorem 1.4}, 
we also obtain
\begin{align}\label{explicit}
\p(n, x, y) 
&=\sum_{\substack{(e_1, \dots, e_n) 
\in \Omega_{x, n}(X) \\ o(e_1)=x, \, t(e_n)=y}}
p(e_1) \cdots p(e_n)\nn \\
&\hspace{1cm} \times \exp\Big( \sum_{i=1}^n
\lambda_*(\bm{x})[ d\Phi_0(\widetilde{e}_i)]_{\Gamma \otimes \mathbb{R}}
\Big)\cdot F_{\bm{x}}\big( \lambda_*(\bm{x})\big)^{-n}\nn\\
&=p(n, x, y) \exp\Big( 
\lambda_*(\bm{x})\big[ \Phi_0(y)-\Phi_0(x)\big]_{\Gamma \otimes \mathbb{R}}
\Big)\cdot F_{\bm{x}}\big( \lambda_*(\bm{x})\big)^{-n}
\end{align}
for every $n \in \mathbb{N}, \, x, y \in V$.
\end{Rem}

\subsection{A relationship to a discrete analogue of Girsanov's theorem}
In this subsection, we discuss a relationship between our formula (\ref{explicit})
and a discrete analogue of Girsanov's theorem due to Fujita \cite{Fujita}.

Let $X=(V, E)$ be a crystal lattice covered with
a one-bouquet graph
$X_0=(V_0, E_0)$; $V_0=\{\bm{x}\}$
and $E_0=\{e, \ol{e}\}$, by the group action
$\Gamma=\la \sigma \ra \cong \mathbb{Z}^1$. 
We consider a random walk on $X_0$ with the transition probability
$$
p(e)=p \quad \text{and} \quad p(\ol{e})=1-p \qquad (0<p<1).
$$
We introduce a bijective linear map 
$\rho_{\mathbb{R}} : \h_1(X_0, \mathbb{R})
\LA \Gamma \otimes \mathbb{R} (\cong \mathbb{R}^1)$
by $\rho_{\mathbb{R}}([e])=\sigma$. 
Then we have 
$\gamma_p=(2p-1)[e]$ and $\rho_{\mathbb{R}}(\gamma_p)=(2p-1)\sigma$. 
Let 
$\{u\} \subset \Hom(\Gamma, \mathbb{R})
=\big(\h^1(X_0, \mathbb{R}), \La \cdot , \cdot \Ra_p\big)$
be a dual basis of $\{\sigma \otimes 1=\sigma\} 
\subset \Gamma \otimes \mathbb{R}$. 
We easily see that $\La u, u\Ra_p=4p(1-p).$ 
Hence the orthogonalization $\{v\} \subset \Hom(\Gamma, \mathbb{R})$ of $\{u\}$
is given by 
$$
v=\frac{1}{\sqrt{4p(1-p)}}u.
$$ 
To the end, we identify $\lambda v \in \Hom(\Gamma, \mathbb{R})$
with $\lambda \in \mathbb{R}$. 
We denote by $\{\vv\} \subset \Gamma \otimes \mathbb{R}$ the dual basis of $\{v\}$.
Then we observe that the realization 
$\Phi_0 : X \LA (\Gamma \otimes \mathbb{R}; \{\vv\})$ defined by
$$
d\Phi_0(\widetilde{e}):=\sigma=\frac{1}{\sqrt{4p(1-p)}}\vv
$$
is the modified standard realization of $X$. 

We now consider the function $F=F_{\bm{x}}(\lambda)$ defined by (\ref{F}), that is, 
$$
F_{\bm{x}}(\lambda)=p\exp\Big(\frac{\lambda}{\sqrt{4p(1-p)}}\Big)
+(1-p)\exp\Big(-\frac{\lambda}{\sqrt{4p(1-p)}}\Big) \qquad (\lambda \in \mathbb{R}).
$$
Then we know that the minimizer $\lambda_*=\lambda_*(\bm{x})$ and $F_{\bm{x}}(\lambda_*)$ are
given  by 
$$
\lambda_* =\sqrt{p(1-p)}\log \frac{p-1}{p}, \qquad F_{\bm{x}}(\lambda_*)=\sqrt{4p(1-p)}.
$$

We fix $x \in V$ satisfying $\Phi_0(x)=\bm{0}$. 
For $y \in V$, we write $\Phi_0(y)=k(y)\vv.$ 
Then the formula (\ref{explicit}) implies
$$
\p(n, x, y)=p(n, x, y) \cdot \Big(\frac{p-1}{p}\Big)^{-k(y)/2} \cdot \big(\sqrt{4p(1-p)}\big)^{-n}
\qquad (n \in \mathbb{N}, \, y \in V).
$$
In  \cite[page 115]{Fujita}, 
the above formula is called
a {\it discrete analogue of 
Girsanov's theorem}  for a non-symmetric random walk 
$\{Z_n\}_{n=0}^\infty$ on $\mathbb{Z}^1$
given by the sum of independent random variables $\{\xi_i\}_{i=1}^\infty$ with
$\mathbb{P}(\xi_i=1)=p$ and $\mathbb{P}(\xi_i=-1)=1-p \, (i=1, 2, \dots)$.
Hence we may regard the formula (\ref{explicit}) as a generalization of the above discrete
Girsanov's theorem to the case of non-symmetric random walks on the $\ell$-bouquet graph.


\section{Examples}

In \cite[Section 7]{IKK}, several examples of the modified standard realization 
of crystal lattices associated with the non-symmetric random walks are discussed.
In this final section, 
we give two concrete examples of non-symmetric random walks on crystal lattices
and calculate the changed transition probability $\p$ for each example.

\subsection{The hexagonal lattice}
\begin{figure}[htbp]
\begin{center}
\includegraphics[width=1\linewidth]{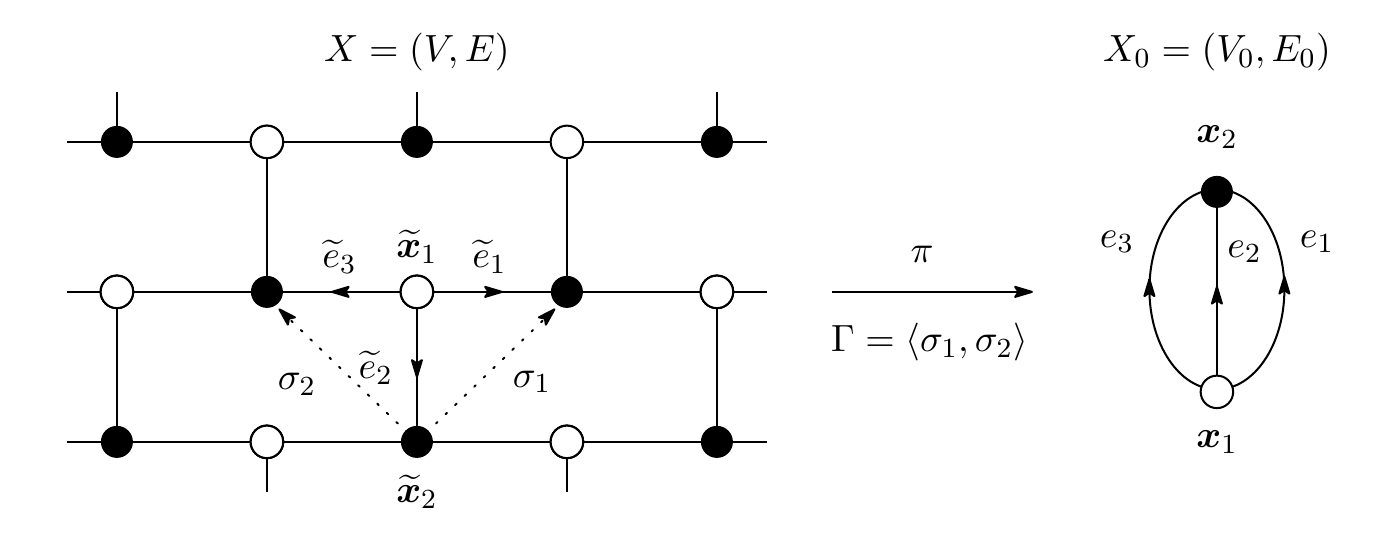}
\caption{Hexagonal lattice and its quotient.}
\label{hexagonal}
\end{center}
\end{figure}

In this subsection, we consider the {\it hexagonal lattice}, 
as a typical example of crystal lattices. 
Let $X=(V, E)$ be a hexagonal lattice, 
$V=\mathbb{Z}^2=\{\x=(x_1, x_2) \, | \, x_1, x_2 \in \mathbb{Z}\}$ and 
$$
E=\big\{(\x, \bm{\mathrm{y}}) \in V^2 \, 
\big| \, \x - \bm{\mathrm{y}}=\pm (1, 0), 
\, \x - \bm{\mathrm{y}}=(0, (-1)^{x_1+x_2})\big\}.
$$
(See Figure \ref{hexagonal}). 
We introduce a non-symmetric random walk on $X$ in the following way.
 If $\x=(x_1, x_2) \in V$ is a vertex so that  $x_1+x_2$ is even, set
$$
p\big(\x, \x+(1, 0)\big)=\frac{1}{2}, 
\quad p\big(\x, \x-(0, 1)\big)=\frac{1}{3}, 
\quad p\big(\x, \x-(1, 0)\big)=\frac{1}{6}.
$$
If $x_1+x_2$ is odd, set
$$
p\big(\x, \x-(1, 0)\big)=\frac{1}{6}, 
\quad p\big(\x, \x+(0, 1)\big)=\frac{1}{3}, 
\quad p\big(\x, \x+(1, 0)\big)=\frac{1}{2}.
$$
We see that $X$ is invariant under the action 
$\Gamma=\la \sigma_1, \sigma_2 \ra ( \cong \mathbb{Z}^2)$ generated by 
$$
\sigma_1(\x)=\x+(1, 1), \quad \sigma_2(\x)=\x+(-1, 1) \qquad (\x \in V).
$$

The quotient graph $X_0=\Gamma \backslash X$
is a finite graph $X_0=(V_0, E_0)$ consisting of two vertices
$\{\bm{x}_1, \bm{x}_2\}$ with three multiple edges 
$E_0=(E_0)_{\bm{x}_1} \cup (E_0)_{\bm{x}_2}
=\{e_1, e_2, e_3\} \cup \{\ol{e}_1, \ol{e}_2, \ol{e}_3\}$. 

We put $[c_1]:=[e_1 * \ol{e}_2]$ and $[c_2]:=[e_3 * \ol{e}_2]$. 
Then, the first homology group 
$\h_1(X_0, \mathbb{R})$ is spanned by $\{[c_1], [c_2]\}$. 
Solving (\ref{invariant}), we have $m(\bm{x}_1)=m(\bm{x}_2)=1/2$. 
We define the surjective linear map 
$\rho_{\mathbb{R}} : \h_1(X_0, \mathbb{R}) 
\LA \Gamma \otimes \mathbb{R} (\cong \mathbb{R}^2)$ by
$\rho_{\mathbb{R}}([c_1]):=\sigma_1, \, \rho_{\mathbb{R}}([c_2]):=\sigma_2.$
Thus, the homological direction $\gamma_p$ 
and the asymptotic direction $\rho_{\mathbb{R}}(\gamma_p)$ are given by
$$
\gamma_p=\frac{1}{6}[c_1]-\frac{1}{6}[c_2], \quad 
\rho_{\mathbb{R}}(\gamma_p)=\frac{1}{6}\sigma_1 - \frac{1}{6}\sigma_2 (\neq \bm{0}),
$$
respectively.
We will determine the modified standard realization 
$\Phi_0 : X \LA (\Gamma \otimes \mathbb{R}, g_0)$. 
We set $\widetilde{\bm{x}}_1=(0, 0), \, \widetilde{\bm{x}}_2=(0, 1)$ in $V$. 
Without loss of generality, 
we may put $\Phi_0(\widetilde{\bm{x}}_1)=\bm{0} \in \Gamma \otimes \mathbb{R}$. 
By (\ref{modified harmonicity}), we have
$$
\Phi_0(\widetilde{\bm{x}}_1)=\bm{0}, 
\quad \Phi_0(\widetilde{\bm{x}}_2)=-\frac{1}{3}\sigma_1-\frac{1}{3}\sigma_2.
$$
Now let $\{v_1, v_2\}$ be an orthonormal basis in 
$\Hom(\Gamma, \mathbb{R})
\big( \subset (\mathcal{H}^1(X_0), \La \cdot , \cdot \Ra_p)\big)$  
and $\{\vv_1, \vv_2\}$ its dual basis in $(\Gamma \otimes \mathbb{R}, g_0)$. Then, we have
$$
\sigma_1=\frac{6\sqrt{7}}{7}\vv_1 - \frac{9\sqrt{70}}{70}\vv_2, 
\quad \sigma_2=\frac{3\sqrt{70}}{10}\vv_2
$$
with the Albanese metric 
$$
\la \sigma_1, \sigma_1 \ra_{g_0}=\frac{63}{10}, \quad 
\la \sigma_1, \sigma_2 \ra_{g_0}=\frac{27}{10}, \quad 
\la \sigma_2, \sigma_2 \ra_{g_0}=\frac{63}{10},
$$
by following the computations in \cite[Subsection 7.3]{IKK}. 
Hence, we find that the modified standard realization 
$\Phi_0 : X \LA (\Gamma \otimes \mathbb{R}, g_0)$ is given by 
$\Phi_0(\widetilde{\bm{x}}_1)=\bm{0}$ and 
$$
\Phi_0(\widetilde{\bm{x}}_2)=-\frac{2\sqrt{7}}{7}\vv_1 - \frac{\sqrt{70}}{7}\vv_2.
$$ 

Now we are in a position to consider the function $F$ introduced in (\ref{F}). 
We identify $\lambda \in \Hom(\Gamma, \mathbb{R})$ 
with $(\lambda_1, \lambda_2) \in \mathbb{R}^2$. Then, we have
\begin{equation}\label{FF}
\begin{split}
F_{\bm{x}_1}(\lambda)&=\frac{1}{2}\exp
\Big(\frac{4\sqrt{7}}{7}\lambda_1 - \frac{19\sqrt{70}}{70}\lambda_2\Big)\\
&\qquad\quad +\frac{1}{3}\exp\Big(-\frac{2\sqrt{7}}{7}\lambda_1 
- \frac{\sqrt{70}}{7}\lambda_2\Big)
+\frac{1}{6}\exp\Big( -\frac{2\sqrt{7}}{7}\lambda_1 
+ \frac{11\sqrt{70}}{70}\lambda_2\Big),\\
F_{\bm{x}_2}(\lambda)&=\frac{1}{6}\exp\Big(-\frac{4\sqrt{7}}{7}\lambda_1 
+ \frac{19\sqrt{70}}{70}\lambda_2\Big)\\
&\qquad\quad +\frac{1}{3}\exp\Big(\frac{2\sqrt{7}}{7}\lambda_1 
+ \frac{\sqrt{70}}{7}\lambda_2\Big)+\frac{1}{2}\exp
\Big( \frac{2\sqrt{7}}{7}\lambda_1 - \frac{11\sqrt{70}}{70}\lambda_2\Big).
\end{split}
\end{equation}
Differentiating both sides of (\ref{FF}) with respect to $\lambda_1$ and $\lambda_2$, we have
$$
\begin{aligned}
\del_1F_{\bm{x}_1}(\lambda)&=\frac{2\sqrt{7}}{7}\exp\Big(\frac{4\sqrt{7}}{7}\lambda_1 
- \frac{19\sqrt{70}}{70}\lambda_2\Big)\\
&\quad -\frac{2\sqrt{7}}{21}\exp\Big(-\frac{2\sqrt{7}}{7}\lambda_1 
- \frac{\sqrt{70}}{7}\lambda_2\Big)-\frac{\sqrt{7}}{21}
\exp\Big( -\frac{2\sqrt{7}}{7}\lambda_1 + \frac{11\sqrt{70}}{70}\lambda_2\Big),\\
\del_2F_{\bm{x}_1}(\lambda)&=-\frac{19\sqrt{70}}{140}
\exp\Big(\frac{4\sqrt{7}}{7}\lambda_1 - \frac{19\sqrt{70}}{70}\lambda_2\Big)\\
&\quad -\frac{\sqrt{70}}{21}\exp\Big(-\frac{2\sqrt{7}}{7}\lambda_1 
- \frac{\sqrt{70}}{7}\lambda_2\Big)+\frac{11\sqrt{70}}{420}
\exp\Big( -\frac{2\sqrt{7}}{7}\lambda_1 + \frac{11\sqrt{70}}{70}\lambda_2\Big),\\
\del_1F_{\bm{x}_2}(\lambda)&=-\frac{2\sqrt{7}}{21}
\exp\Big(-\frac{4\sqrt{7}}{7}\lambda_1 + \frac{19\sqrt{70}}{70}\lambda_2\Big)\\
&\quad +\frac{2\sqrt{7}}{21}\exp\Big(\frac{2\sqrt{7}}{7}\lambda_1 
+ \frac{\sqrt{70}}{7}\lambda_2\Big)+\frac{\sqrt{7}}{7}
\exp\Big( \frac{2\sqrt{7}}{7}\lambda_1 - \frac{11\sqrt{70}}{70}\lambda_2\Big),\\
\del_2F_{\bm{x}_2}(\lambda)&=\frac{19\sqrt{70}}{420}
\exp\Big(-\frac{4\sqrt{7}}{7}\lambda_1 + \frac{19\sqrt{70}}{70}\lambda_2\Big)\\
&\quad +\frac{\sqrt{70}}{21}\exp\Big(\frac{2\sqrt{7}}{7}\lambda_1 
+ \frac{\sqrt{70}}{7}\lambda_2\Big)-\frac{11\sqrt{70}}{140}
\exp\Big( \frac{2\sqrt{7}}{7}\lambda_1 - \frac{11\sqrt{70}}{70}\lambda_2\Big).
\end{aligned}
$$
To find the minimizers of functions $F_{\bm{x}_1}(\cdot)$ and $F_{\bm{x}_2}(\cdot)$, 
it is sufficient to solve the following two algebraic equations:
$$
\begin{cases}
\del_1F_{\bm{x}_1}(\lambda_1, \lambda_2)=0 \\
\del_2F_{\bm{x}_1}(\lambda_1, \lambda_2)=0 
\end{cases}, \quad \begin{cases}
\del_1F_{\bm{x}_2}(\lambda_1, \lambda_2)=0 \\
\del_2F_{\bm{x}_2}(\lambda_1, \lambda_2)=0 
\end{cases}.
$$
Solving these equations, the minimizers $\lambda_*(\bm{x}_1)$, $\lambda_*(\bm{x}_2)$
of $F_{\bm{x}_1}(\cdot)$, $F_{\bm{x}_2}(\cdot)$ are given by
\begin{equation}\label{minimizer}
\begin{split}
\lambda_*(\bm{x}_1)&=\Big(\frac{\sqrt{7}}{14}\log 26 + \frac{\sqrt{7}}{6}\log \frac{14}{3}, 
\frac{\sqrt{70}}{21}\log 26\Big),\\
 \lambda_*(\bm{x}_2)&=\Big( \frac{\sqrt{7}}{14}\log\frac{3}{26} - \frac{\sqrt{7}}{6}\log 14, 
 \frac{\sqrt{70}}{21}\log\frac{3}{26}\Big),
 \end{split}
 \end{equation}
respectively. Then, it follows from (\ref{FF}) and (\ref{minimizer}) that
$$
F_{\bm{x}_1}\big(\lambda_*(\bm{x}_1)\big)
=7 \cdot 26^{-13/21} \cdot \Big(\frac{14}{3}\Big)^{-1/3},
\quad
F_{\bm{x}_2}\big(\lambda_*(\bm{x}_2)\big)
=\frac{21}{26} \cdot 14^{-1/3} \cdot \Big(\frac{3}{26}\Big)^{-8/21}.
$$
Finally, we determine the changed transition probability $\p$ by 
$$
\p(e_1)=\frac{1}{3}, \quad \p(e_2)=\frac{1}{21}, \quad \p(e_3)=\frac{13}{21}, \quad 
\p(\ol{e}_1)=\frac{1}{3}, \quad \p(\ol{e}_2)=\frac{1}{21}, \quad \p(\ol{e}_3)=\frac{13}{21}.
$$
Then, the invariant measure $\m : \{\bm{x}_1, \bm{x}_2\} \LA (0, 1]$ 
is also given by $\m(\bm{x}_1)=\m(\bm{x}_2)=1/2$. 
Hence we know that the random walk 
associated with the changed transition probability $\p$ is $\m$-symmetric, that is,
$$
\p(e)\m\big(o(e)\big)=\p(\ol{e})\m\big(t(e)\big) \qquad (e \in E_0).
$$
It automatically implies $\gamma_{\p}=0$ and also $\rho_{\mathbb{R}}(\gamma_{\p})=\bm{0}$.

\subsection{The dice lattice}

\begin{figure}[htbp]
\begin{center}
\includegraphics[width=1\linewidth]{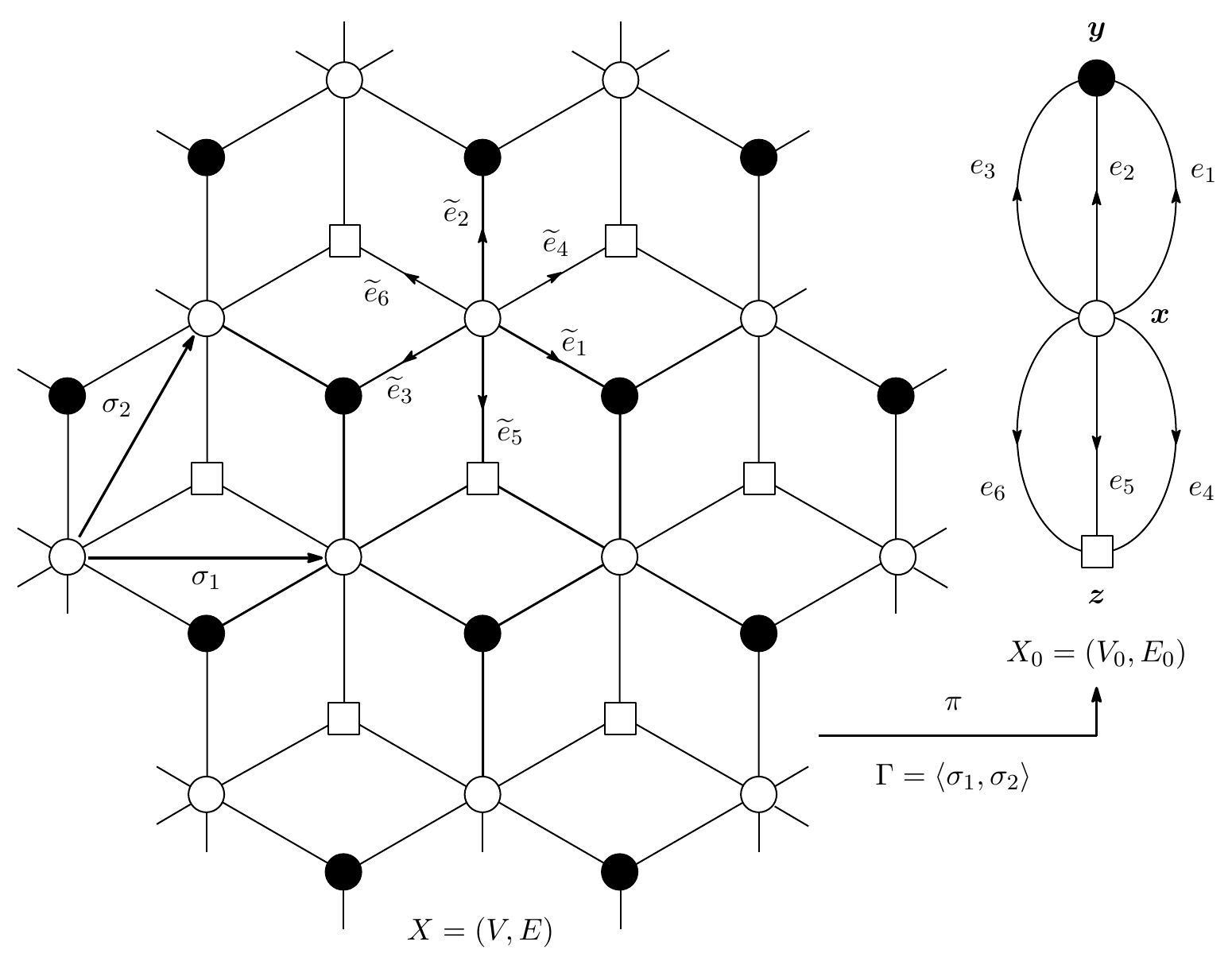}
\caption{Dice lattice and its quotient.}
\label{Dice}
\end{center}
\end{figure}

In this subsection, we discuss a non-symmetric random walk on an infinite graph 
called the {\it dice lattice} or the {\it dice graph}. 
The dice graph $X=(V, E)$ is one of abelian covering graphs 
which has a free action by the lattice group  $\Gamma \cong \mathbb{Z}^2$ 
generated by $\sigma_1, \sigma_2$, and 
the corresponding quotient graph $X_0=(V_0, E_0):=\Gamma \backslash X$ 
is a finite graph consisting of three vertices $V_0=\{\bm{x}, \bm{y}, \bm{z}\}$. 
(See Figure \ref{Dice} or the description in \cite{Sunada}.) 
In view of the shape of the quotient graph $X_0$, 
we may regard the dice lattice $X$ as something
like a ``{\it hybrid\,}'' of a triangular lattice and a hexagonal lattice. 

From now on, we consider a non-symmetric random walk on $X$ 
by giving the transition probability on the quotient $X_0$
 in the following way. We set
\begin{align}
p(e_1)&=\frac{1}{4}, & p(e_2)&=\frac{1}{6}, & p(e_3)&=\frac{1}{12}, & p(e_4)&=\frac{1}{4}, & 
p(e_5)&=\frac{1}{6}, & p(e_6)&=\frac{1}{12},\nn\\
p(\ol{e}_1)&=\frac{1}{6}, & p(\ol{e}_2)&=\frac{1}{3}, & p(\ol{e}_3)&=\frac{1}{2}, & p(\ol{e}_4)&=\frac{1}{6}, &
p(\ol{e}_5)&=\frac{1}{3}, & p(\ol{e}_6)&=\frac{1}{2}.\nn
\end{align}
Solving (\ref{invariant}), we have $m(\bm{x})=1/2, \, m(\bm{y})=m(\bm{z})=1/4$. 

Next we define four 1-cycles $[c_1], [c_2], [c_3], [c_4]$ on $X_0$ by 
$$
[c_1]:=[e_1 * \ol{e}_2], \quad [c_2]:=[e_3 * \ol{e}_2], \quad [c_3]:=[e_4 * \ol{e}_5], \quad [c_4]:=[e_6 * \ol{e}_5].
$$
Then $\{[c_1], [c_2], [c_3], [c_4]\}$ spans
 the first homology group $\h_1(X_0, \mathbb{R})$.
We define the linear map $\rho_{\mathbb{R}}$ 
from $\h_1(X_0, \mathbb{R})$ onto $\Gamma \otimes \mathbb{R} \cong \mathbb{R}^2$ by
$$
\rho_{\mathbb{R}}([c_1])=\sigma_1 -\sigma_2, \quad 
\rho_{\mathbb{R}}([c_2])=-\sigma_2, \quad 
\rho_{\mathbb{R}}([c_3])=\sigma_2, \quad 
\rho_{\mathbb{R}}([c_4])=\sigma_2 - \sigma_1.
$$
Then, the homological direction $\gamma_p$ and the asymptotic direction 
$\rho_{\mathbb{R}}(\gamma_p)$
of the random walk on $X_0$ 
are calculated as
$$
\begin{aligned}
\gamma_p&=\frac{1}{12}\big( [c_1] - [c_2] + [c_3] - [c_4] \big),\\
\rho_{\mathbb{R}}(\gamma_p)&=
\frac{1}{12}\big\{ (\sigma_1 - \sigma_2)-(-\sigma_2)+\sigma_2
- (\sigma_2 - \sigma_1)\big\}=\frac{1}{6}\sigma_1 (\neq \bm{0}),
\end{aligned}
$$
respectively. 

Now we determine the modified standard realization $\Phi_0 : X \LA \Gamma \otimes \mathbb{R}$. 
Here we may put $\Phi_0\big(o(\widetilde{e}_i)\big)=\bm{0} \,\, (i=1, 2, 3, 4, 5, 6)$,
 without loss of generality.
Noting (\ref{modified harmonicity}) and the group action $\Gamma$, we have
\begin{align}\label{phi}
\Phi_0\big(t(\widetilde{e}_1)\big)&=\frac{2}{3}\sigma_1 - \frac{1}{3}\sigma_2, 
& \Phi_0\big(t(\widetilde{e}_2)\big)&=-\frac{1}{3}\sigma_1+\frac{2}{3}\sigma_2, \nn\\
\Phi_0\big(t(\widetilde{e}_3)\big)&=-\frac{1}{3}\sigma_1 - \frac{1}{3}\sigma_2, 
& \Phi_0\big(t(\widetilde{e}_4)\big)&=\frac{1}{3}\sigma_1+\frac{1}{3}\sigma_2, \\
\Phi_0\big(t(\widetilde{e}_5)\big)&=\frac{1}{3}\sigma_1 - \frac{2}{3}\sigma_2, 
& \Phi_0\big(t(\widetilde{e}_6)\big)&=-\frac{2}{3}\sigma_1+\frac{1}{3}\sigma_2.\nn
\end{align}

Let $\{\omega_1, \omega_2, \omega_3, \omega_4\} 
\subset \h^1(X_0, \mathbb{R})$ be the dual basis of 
$\{[c_1], [c_2], [c_3], [c_4]\}$, 
that is, $\omega_i([c_j])=\delta_{ij} \, (1 \leq i, j \leq 4)$. 
Recalling that each $\omega_i$ is a modified harmonic 1-form, we have
$$
\omega_1(e_1)=\frac{3}{4},    
\quad  \omega_1(e_2)=\omega_1(e_3)=-\frac{1}{4}, 
\quad
\omega_1(e_4)=\omega_1(e_5)=\omega_1(e_6)=-\frac{1}{12}, 
$$
$$
\omega_2(e_1)=\omega_2(e_2)=-\frac{5}{12},    
\quad  \omega_2(e_3)=\frac{7}{12}, 
\quad
\omega_2(e_4)=\omega_2(e_5)=\omega_2(e_6)=\frac{1}{12}, 
$$
$$
\omega_3(e_1)=\omega_3(e_2)=\omega_3(e_3)=-\frac{1}{12},    
\quad  \omega_3(e_4)=\frac{3}{4}, 
\quad
\omega_3(e_5)=\omega_3(e_6)=-\frac{1}{4}, 
$$
$$
\omega_4(e_1)=\omega_4(e_2)=\omega_4(e_3)=\frac{1}{12},    
\quad  \omega_4(e_4)=\omega_4(e_5)=-\frac{5}{12}, 
\quad
\omega_4(e_6)=\frac{7}{12}.
$$
Then, the direct computation gives us
\begin{align}
\La \omega_1, \omega_1 \Ra&=\frac{1}{9}, 
& \La \omega_1, \omega_2 \Ra &=-\frac{1}{18}, 
& \La \omega_1, \omega_3 \Ra&=-\frac{1}{72},
&\La \omega_1, \omega_4 \Ra&=\frac{1}{72}, 
& \La \omega_2, \omega_2 \Ra &= \frac{1}{9}, \nn\\
\La \omega_2, \omega_3 \Ra&=\frac{1}{72}, 
& \La \omega_2, \omega_4 \Ra &=-\frac{1}{72}, 
&\La \omega_3, \omega_3 \Ra&=\frac{1}{9},
& \La \omega_3, \omega_4 \Ra&=-\frac{1}{18}, 
& \La \omega_4, \omega_4 \Ra &= \frac{1}{9}.\nn
\end{align}
We introduce the basis $\{u_1, u_2\}$ in $\Hom(\Gamma, \mathbb{R})$ 
by the dual of $\{\sigma_1 \otimes 1, \sigma_2 \otimes 1\} 
\subset \Gamma \otimes \mathbb{R}$. 
Since the dice graph $X$ is a non-maximal 
abelian covering graph of $X_0$ with 
$\Gamma \cong \la \sigma_1, \sigma_2 \ra$,
we need to find a $\mathbb{Z}$-basis of the lattice
$$
L=\big\{ \omega \in \h^1(X_0, \mathbb{R}) \, 
\big| \, \omega([c])=0 \text{ for every cycle $\widetilde{c}$ on $X$} \big\}.
$$
It is easy to find that 
$$
u_1={}^t \rho_{\mathbb{R}}(u_1)=\omega_1 - \omega_4, \quad 
u_2={}^t \rho_{\mathbb{R}}(u_2)=-\omega_1 - \omega_2 + \omega_3 + \omega_4
$$
form a $\mathbb{Z}$-basis of the lattice $L$. 
Carrying out the direct computation again, we have
\begin{equation}\label{u}
\La u_1, u_1 \Ra = \frac{7}{36}, 
\quad \La u_1, u_2 \Ra = -\frac{1}{9}, 
\quad \La u_2, u_2 \Ra =\frac{2}{9}.
\end{equation}
Let $\{v_1, v_2\}$ be the Gram--Schmidt orthogonalization 
of the basis $\{u_1, u_2\} \subset \Hom(\Gamma, \mathbb{R})$. 
By (\ref{u}), we have
$$
v_1=\frac{6\sqrt{7}}{7}u_1, \quad 
v_2=\frac{6\sqrt{70}}{35}u_1 + \frac{3\sqrt{70}}{10}u_2.
$$
We denote by $\{\vv_1, \vv_2\}$ 
the dual basis of $\{v_1, v_2\}$ in $\Gamma \otimes \mathbb{R}$. 
Then we obtain
\begin{equation}\label{sigma}
\sigma_1=\frac{6\sqrt{7}}{7}\vv_1, \quad \sigma_2
=\frac{6\sqrt{70}}{35}\vv_1 + \frac{3\sqrt{70}}{10}\vv_2
\end{equation}
Combining (\ref{sigma}) with (\ref{phi}),
we finally determine the modified standard realization 
$\Phi_0 : X \LA (\Gamma \otimes \mathbb{R}, g_0) \cong (\mathbb{R}^2; \{\vv_1, \vv_2\})$ 
of $X$ by
\begin{align}
\Phi_0\big(t(\widetilde{e}_1)\big) &= \frac{20\sqrt{7}-2\sqrt{70}}{35}\vv_1 - \frac{\sqrt{70}}{10}\vv_2, & 
\Phi_0\big(t(\widetilde{e}_2)\big) &= \frac{-10\sqrt{7}+4\sqrt{70}}{35}\vv_1 + \frac{\sqrt{70}}{5}\vv_2, \nn\\
\Phi_0\big(t(\widetilde{e}_3)\big) &= \frac{-10\sqrt{7}-2\sqrt{70}}{35}\vv_1 - \frac{\sqrt{70}}{10}\vv_2, & 
\Phi_0\big(t(\widetilde{e}_4)\big) &= \frac{10\sqrt{7}+2\sqrt{70}}{35}\vv_1 + \frac{\sqrt{70}}{10}\vv_2, \nn\\
\Phi_0\big(t(\widetilde{e}_5)\big) &= \frac{10\sqrt{7}-4\sqrt{70}}{35}\vv_1 - \frac{\sqrt{70}}{5}\vv_2, & 
\Phi_0\big(t(\widetilde{e}_6)\big) &= \frac{-20\sqrt{7}+2\sqrt{70}}{35}\vv_1 + \frac{\sqrt{70}}{10}\vv_2. \nn
\end{align}
(See Figure \ref{dice-standard} below.)

\begin{figure}[htbp]
\begin{center}
\includegraphics[width=1\linewidth]{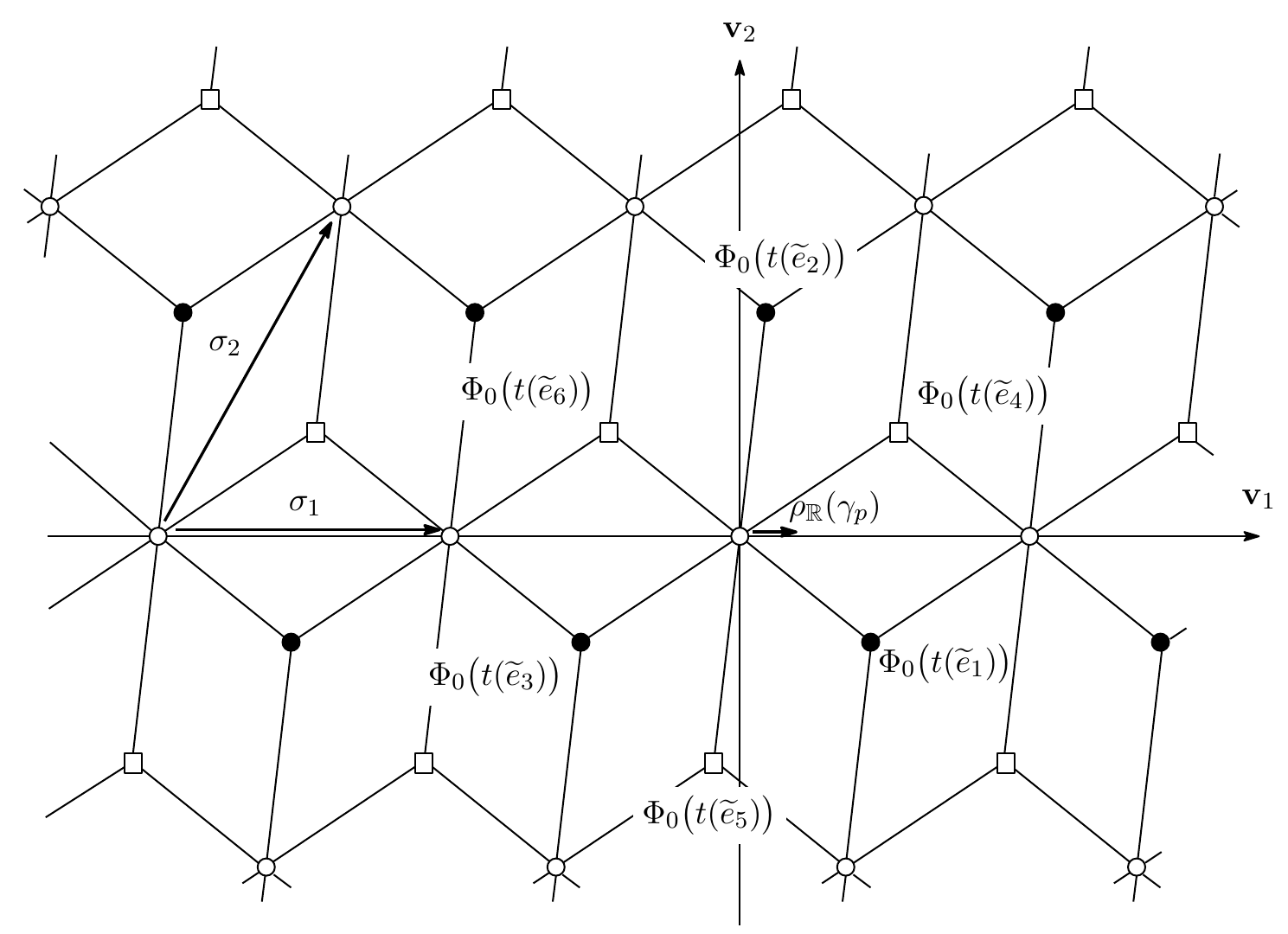}
\caption{Modified standard realization of the dice lattice.}
\label{dice-standard}
\end{center}
\end{figure}

Now we are in a position to consider the function 
$F : V_0 \times \Hom(\Gamma, \mathbb{R}) \LA (0, \infty)$ 
defined by (\ref{F}).  In what follows, 
we identify $\lambda=\lambda_1 v_1 + \lambda_2 v_2 \in \Hom(\Gamma, \mathbb{R})$ 
with $(\lambda_1, \lambda_2) \in \mathbb{R}^2$. 
In addition, it follows from $\Phi_0\big(o(\widetilde{e}_i)\big)=\bm{0} \, (i=1, 2, 3, 4, 5 ,6)$ that 
\begin{align}\label{dphi2}
d\Phi_0(\widetilde{e}_1) &= \frac{20\sqrt{7}-2\sqrt{70}}{35}\vv_1 - \frac{\sqrt{70}}{10}\vv_2, & 
d\Phi_0(\widetilde{e}_2) &= \frac{-10\sqrt{7}+4\sqrt{70}}{35}\vv_1 + \frac{\sqrt{70}}{5}\vv_2, \nn\\
d\Phi_0(\widetilde{e}_3) &= \frac{-10\sqrt{7}-2\sqrt{70}}{35}\vv_1 - \frac{\sqrt{70}}{10}\vv_2, & 
d\Phi_0(\widetilde{e}_4) &= \frac{10\sqrt{7}+2\sqrt{70}}{35}\vv_1 + \frac{\sqrt{70}}{10}\vv_2, \\
d\Phi_0(\widetilde{e}_5) &= \frac{10\sqrt{7}-4\sqrt{70}}{35}\vv_1 - \frac{\sqrt{70}}{5}\vv_2, & 
d\Phi_0(\widetilde{e}_6) &= \frac{-20\sqrt{7}+2\sqrt{70}}{35}\vv_1 + \frac{\sqrt{70}}{10}\vv_2 \nn
\end{align}
Then, by (\ref{dphi2}), we have
$$
\begin{aligned}
&F_{\bm{x}}(\lambda_1, \lambda_2) \\
&=
   \frac{1}{4}\exp\Big(\frac{20\sqrt{7}-2\sqrt{70}}{35}\lambda_1 - \frac{\sqrt{70}}{10}\lambda_2\Big)
   +\frac{1}{6}\exp\Big( \frac{-10\sqrt{7}+4\sqrt{70}}{35}\lambda_1 + \frac{\sqrt{70}}{5}\lambda_2\Big)\\
   &\quad+\frac{1}{12}\exp\Big( \frac{-10\sqrt{7}-2\sqrt{70}}{35}\lambda_1 - \frac{\sqrt{70}}{10}\lambda_2\Big)
   +\frac{1}{4}\exp\Big( \frac{10\sqrt{7}+2\sqrt{70}}{35}\lambda_1 + \frac{\sqrt{70}}{10}\lambda_2\Big)\\
   &\quad+\frac{1}{6}\exp\Big( \frac{10\sqrt{7}-4\sqrt{70}}{35}\lambda_1 - \frac{\sqrt{70}}{5}\lambda_2\Big)
   +\frac{1}{12}\exp\Big( \frac{-20\sqrt{7}+2\sqrt{70}}{35}\lambda_1 + \frac{\sqrt{70}}{10}\lambda_2\Big),\\
&F_{\bm{y}}(\lambda_1, \lambda_2)\\
&=\frac{1}{6}\exp\Big( \frac{-20\sqrt{7}+2\sqrt{70}}{35}\lambda_1 + \frac{\sqrt{70}}{10}\lambda_2\Big)
+\frac{1}{3}\exp\Big( \frac{10\sqrt{7}-4\sqrt{70}}{35}\lambda_1 - \frac{\sqrt{70}}{5}\lambda_2\Big)\\
&\quad +\frac{1}{2}\exp\Big( \frac{10\sqrt{7}+2\sqrt{70}}{35}\lambda_1 + \frac{\sqrt{70}}{10}\lambda_2\Big),\\
&F_{\bm{z}}(\lambda_1, \lambda_2)\\
&=\frac{1}{6}\exp\Big( \frac{-10\sqrt{7}-2\sqrt{70}}{35}\lambda_1 - \frac{\sqrt{70}}{10}\lambda_2\Big)
+\frac{1}{3}\exp\Big( \frac{-10\sqrt{7}+4\sqrt{70}}{35}\lambda_1 + \frac{\sqrt{70}}{5}\lambda_2\Big)\\
&\quad +\frac{1}{2}\exp\Big( \frac{20\sqrt{7}-2\sqrt{70}}{35}\lambda_1 - \frac{\sqrt{70}}{10}\lambda_2\Big).
\end{aligned}
$$
By solving the following equations:
$$
\begin{cases}
\del_1F_{\bm{x}}(\lambda_1, \lambda_2)=0 \\
\del_2F_{\bm{x}}(\lambda_1, \lambda_2)=0 
\end{cases}, \quad \begin{cases}
\del_1F_{\bm{y}}(\lambda_1, \lambda_2)=0 \\
\del_2F_{\bm{y}}(\lambda_1, \lambda_2)=0 
\end{cases}, \quad \begin{cases}
\del_1F_{\bm{z}}(\lambda_1, \lambda_2)=0 \\
\del_2F_{\bm{z}}(\lambda_1, \lambda_2)=0 
\end{cases},
$$
we find that  the minimizers $\lambda_*(\bm{x}), \lambda_*(\bm{y})$ and $\lambda_*(\bm{z})$
of functions $F_{\bm{x}}(\cdot), F_{\bm{y}}(\cdot)$ and $F_{\bm{z}}(\cdot)$ are given by
$$
\begin{aligned}
\lambda_*(\bm{x}) 
&= \Big( -\frac{\sqrt{7}}{6}\log 3, \frac{4\sqrt{7}-\sqrt{70}}{42}\log 3 \Big),\\
\lambda_*(\bm{y}) 
&= \Big( -\frac{\sqrt{7}}{6}\log 3, \frac{\sqrt{70}}{21}\log 2 + \frac{2\sqrt{7}-\sqrt{70}}{21}\log 3\Big),\\
\lambda_*(\bm{z}) 
&= \Big( -\frac{\sqrt{7}}{6}\log 3, -\frac{\sqrt{70}}{21}\log 2 + \frac{2\sqrt{7}}{21}\log 3\Big),
\end{aligned}
$$
respectively, and 
$$
F_{\bm{x}}\big(\lambda_*(\bm{x})\big)=\frac{\sqrt{3}+1}{3}, \quad 
F_{\bm{y}}\big(\lambda_*(\bm{y})\big)=F_{\bm{z}}\big(\lambda_*(\bm{z})\big)=3 \cdot 6^{-2/3}.
$$

Finally, we determine the changed transition probability $\p$ by 
$$
\p(e_1)=\frac{3-\sqrt{3}}{8}, \quad \p(e_2)=\frac{\sqrt{3}-1}{4}, \quad \p(e_3)=\frac{3-\sqrt{3}}{8}, 
$$
$$
 \p(e_4)=\frac{3-\sqrt{3}}{8}, \quad 
\p(e_5)=\frac{\sqrt{3}-1}{4}, \quad \p(e_6)=\frac{3-\sqrt{3}}{8},
$$
$$
\p(\ol{e}_1)=\p(\ol{e}_2)=\p(\ol{e}_3)=\p(\ol{e}_4)=\p(\ol{e}_5)=\p(\ol{e}_6)=\frac{1}{3}.
$$
Then, the invariant measure $\m$ is given by $\m(\bm{x})=1/2$ and $\m(\bm{y})=\m(\bm{z})=1/4$. 
Moreover, the homological direction $\gamma_\p$ and
the asymptotic direction $\rho_{\mathbb{R}}(\gamma_\p)$
 of the random walk associated with $\p$ are given by
$$
\gamma_{\p}=\frac{5-3\sqrt{3}}{48}\big([c_1]+[c_2]+[c_3]+[c_4]\big) (\neq 0),
\quad \rho_{\mathbb{R}}(\gamma_\p)=\bm{0},
$$
respectively. Then it follows from Theorem \ref{asymptotic-p} that 
there exist some positive constants $C_1$ and $C_2$ such that
$$
C_1p(n, x, y)\big( \sqrt[6]{12}(\sqrt{3}-1)\big)^n \leq \p(n, x, y)
\leq C_2p(n, x, y)\big( \sqrt[6]{12}(\sqrt{3}-1)\big)^n 
$$
for all $n \in \mathbb{N}$ and $x, y \in V$. 

\vspace{5mm}
\noindent
{\bf Acknowledgement.} The author would like to thank his advisor 
Professor Hiroshi Kawabi for useful discussion and constant encouragement.
He expresses his gratitude to Professors Satoshi Ishiwata, Tomoyuki Kakehi, Atsushi Katsuda
and Ryokichi Tanaka
for valuable advice and comments. 
He also would like to thank the anonymous referee for providing helpful comments and suggestions.  

\end{document}